\documentclass[12pt, A4paper]{article}

\usepackage[utf8]{inputenc}
\usepackage{amsthm,amsmath,amssymb,graphicx,enumerate,comment}
\usepackage{tikz}
\usepackage{xcolor,colortbl}

\newtheorem{theorem}{Theorem}[section]
\newtheorem{lemma}[theorem]{Lemma}
\newtheorem{proposition}[theorem]{Proposition}

\newtheorem{conjecture}{Conjecture}

\theoremstyle{definition}

\newtheorem{assumption}[theorem]{Assumption}
\newtheorem{example}[theorem]{Example}
\newtheorem{remark}[theorem]{Remark}

\newcommand{\ZZ}{\mathbb{Z}}

\newcommand{\arxiv}[1]{\texttt{arXiv:#1}} 

\title{All generalized rose window graphs are hamiltonian }
\author{Simona\,Bonvicini\thanks{Dipartimento di Scienze Fisiche, Informatiche e Matematiche, 
Universit\`a di Modena e Reggio Emilia, via Campi 213/b, 41126 Modena, Italy.},
\ Toma\v z\, Pisanski\thanks{University of Primorska  and Institute of Mathematics, Physics and Mechanics, Slovenia},
\ Arjana\,\v Zitnik\thanks{University of Ljubljana and Institute of Mathematics, Physics and Mechanics, Ljubljana, Slovenia}}

\date{\today}

\begin{document}

\maketitle
\centerline{Dedicated to Brian Alspach, who opened this research area.} 

\begin{abstract}
A \emph{bicirculant} is a regular, $d$-valent graph that admits a semiregular automorphism of order $m$ having two vertex-orbits of size $m$. 
The vertices of each orbit induce a circulant graph of order $m$ and the remaining edges span a regular bipartite graph of valence, say $s$, $1 \leq s \leq d$, connecting the two vertex-orbits. Generalized Petersen graphs constitute a prominent family of bicirculants, with $d = 3$ and $s = 1$. 
In 1983, Brian Alspach proved that all generalized Petersen graphs are hamiltonian, except for the family $G(m, 2)$ with $m\equiv 5\pmod 6$.
In this paper we conjecture that among all connected bicirculants of valence at least 2, there are no other exceptions.
It follows from various sources that the conjecture is true for all cubic bicirculants. 
In this paper we prove the conjecture for quartic bicirulants with $s = 2$, also known as the generalized rose window graphs.
\medskip

\medskip

\noindent
\textbf{Math Subj Class (2020)}: 
05C45, 
05C25, 
05C76,  
05C70,  
05E18. 
\end{abstract} 

\section{Introduction}
Motivated by the inspiring seminal work of Brian Alspach on generalized Petersen graphs \cite{Al11983} and the subsequent papers on the hamiltonian properties of certain families of cubic graphs  \cite{AlspachLiu,AlZh1989}, we address the problem of existence of a Hamilton cycle in a larger class of bicirculant graphs.

A \emph{bicirculant} is a regular, $d$-valent graph that admits a semiregular automorphism of order $m$ having two vertex-orbits of size $m$. The vertices of each orbit induce a circulant graph and the remaining edges span a regular bipartite graph of valence, say $s, 1 \leq s \leq d$, connecting the two orbits.  Formal definitions are given in Section  \ref{sec:properties}. 

In general, a regular graph that admits a semiregular automorphism with $k \ge 1$  vertex-orbits is called a \emph{polycirculant} or sometimes a \emph{multicirculant} \cite{BPpolyciclic,kmulticirculant}. Polycirculants with $k = 1$ vertex-orbits are the circulants. 
Bicirculants therefore constitute the next case where $k = 2$. While it is relatively easy to show that all circulants are hamiltonian, see \cite{MarusicCayleyAbelian}, the problem which bicirculants are hamiltonian is still widely open. In particular, it is not even known whether all Cayley graphs on dihedral groups are hamiltonian \cite{AlChDe2010}. Note that all Cayley graphs on dihedral groups are bicirculants.

The generalized Petersen graphs are clearly bicirculants. 
The rims determine the orbits and the spokes constitute a matching between them. 
In \cite{Al11983} Brian Alspach  classified the hamiltonian generalized Petersen graphs: he proved that among the generalized Petersen graphs only the graphs $G(m, 2)$ with $m\equiv 5\pmod 6$  are not  hamiltonian. 
In this paper we  pose the following conjecture: 

\begin{conjecture} \label{conjectB}
Every connected bicirculant, except for the $K_2$ and the generalized Petersen graphs $G(m, 2)$ with $m\equiv 5\pmod 6$, is hamiltonian. 
\end{conjecture}

As bicirculants, generalized Petersen graphs have parameters $d = 3$ and $s = 1$. However, the whole class of bicirculants with parameters $d = 3, s = 1$ consists of $I$-graphs, first introduced in the Foster Census \cite{Foster}. The classification of hamiltonian generalized Petersen graphs from \cite{Al11983}  was extended to $I$-graphs in 2017 \cite{BoPi2017}. It has been proven that all proper $I$-graphs are hamiltonian. The conjecture therefore holds for all bicirculants with parameters $d = 3, s = 1$. Cubic bicirculants fall into three classes, depending on $s$, with $s=1,2,3$ \cite{Pi2007}. Alspach and Zhang \cite{AlZh1989} dealt with the case $d=s=3$. Note that bicirculants with $d = s$ are known as cyclic Haar graphs; they are a special class of Cayley graphs on dihedral groups \cite{HlMaPi2002}.  This essentially covered all connected cubic bicirculants.

In this paper, we take the next step in attacking the case $d = 4$ by resolving the subcase $s = 2$. 
The bicirculants with parameters $d=4$ and $s=2$ are called \emph{generalized rose window graphs} \cite{DaBaPiZi2025}. 
The following is our main result.

\begin{theorem}  \label{th_rosewin_hamiltonian}
Every connected generalized rose window graph is hamiltonian.
\end{theorem}

The rose window graphs, which were introduced by Steve Wilson in 2008 \cite{wilsonRW}, are contained in the family of generalized rose window graphs. 
Informally, a rose window graph is obtained from a generalized Petersen graph by adding an additional set of spokes to its edge set that preserves the semiregular symmetry. Rose window graphs turned out to be a very interesting family of graphs. As they belong to the class of bicirculant graphs, they have many symmetries. Some are vertex-transitive \cite{RWvt} and some are even edge-transitive \cite{RWet} or Cayley \cite{RWCayley}. 
In addition, several of their other properties were studied: isomorphsms \cite{RWiso}, domination \cite{RWdomination}, relation to maps \cite{RWatmaps,RWrotmaps}, etc.

The relationship between generalized rose window graphs and rose window graphs is analogous to the relationship between $I$-graphs and generalized Petersen graphs. 
While generalized Petersen graphs and rose window graphs are necessarily connected, the $I$-graphs and generalized rose window graphs need not be connected. Moreover, the removal of a matching consisting of spokes from a connected generalized rose window graph   results in a disconnected graph whose connected components are $I$-graphs, unless certain arithmetic conditions that will be specified in the next section are satisfied. 
Because of this, the existence of a Hamilton cycle in a generalized rose window graph cannot easily follow from the existence of a Hamilton cycle in 
$I$-graphs. To prove that all generalized rose window graphs are hamiltonian we had to develop several completely novel tools, which are potentially useful for constructing Hamilton cycles in larger families of bicirculants.

From the point of symmetry,  the analogy is more intricate.
Although some proper $I$-graphs may possess symetries  not present in any generalized Petersen graph, none of them is vertex-transitive  \cite{BPZ2005}. 
However, there exist generalized rose window graphs that are Cayley graphs and others that are vertex-transitive and non-Cayley \cite{DaBaPiZi2025}.

Therefore, the results presented in this paper are also important in the context of the Lov\'asz conjecture, a variant of which can be stated as: Every finite connected vertex-transitive graph, except for the five known exceptions, is hamiltonian \cite{MarusicCayleyAbelian}.
By Theorem \ref{th_rosewin_hamiltonian} we confirm this conjecture within the class of vertex-transitive generalized rose window graphs. 

The paper is organized as follows. In Section \ref{sec:properties} we give a formal definition of bicirculants and review the basic properties of bicirculants, and in particular generalized rose window graphs. We recall the notions of generalized Petersen graphs and $I$-graphs \cite{Foster,BPZ2005}, as these graphs appear as subgraphs, actually as spanning sub-bicirculants, of the generalized rose window graphs. In Section \ref{sec:Igraphs}  we also review some known results on hamiltonicity of these graphs. 
We classify  Hamilton cycles in I-graphs into three types. Each of the three types is then used in a different construction in Section \ref{sec_rosewin}, where we prove our main result that all generalized rose window graphs are hamiltonian.

The proof can be briefly described as follows. By removing a suitable matching from a generalized rose window graph $G$ we obtain an $I$-graph $H$, which can be connected or disconnected. As shown in Section \ref{sec:Igraphs}, every connected component of this $I$-graph $H$ contains a Hamilton  cycle or path of a special type. These structures provide subpaths that can be combined into a Hamilton cycle of the entire graph $G$ by using some of the removed edges.

In the last section we then discuss the hamiltonian problem for more general bicirculant graphs.  As a consequence of Theorem \ref{th_rosewin_hamiltonian}, combined with the results from \cite{AlZh1989}, we obtain that every connected bicirculant with $d \ge 5$ and $s = d - 2$ is hamiltonian if $m$ is a product of at most  three prime powers. In particular, this is true for the Taba\v cjn graphs \cite{Taba1,Taba2},  pentavalent bicirculants with $s=3$.

\section{Bicirculants and their properties}
\label{sec:properties}

In this section we give a formal definition of bicirculants, generalized rose window graphs and $I$-graphs, and recall some of their properties.

A \emph{bicirculant} can be described as follows. Given an integer $m \ge 1$ and sets $R,S,T \subseteq \ZZ_m$ such that $R=-R$, $T=-T$, 
$0 \not\in R \cup T$ and $0 \in S$, the graph $B(m;R,S,T)$ has vertex set 
$V=V_1 \cup V_2$, where $V_1=\{u_0,\dots,u_{m-1}\}$ and $V_2=\{v_0,\dots,v_{m-1}\}$, and edge set 
$$E=\{u_iu_{i+j}| \ i \in\ZZ_m, j \in R\} \cup 
      \{v_iv_{i+j}| \ i \in\ZZ_m, j \in T\} \cup
      \{u_iv_{i+j}| \ i \in\ZZ_m, j \in S\}.$$ 
Obviously, the mapping $\alpha:V \to V$, defined by $\alpha(u_i)=u_{i+1}$, $\alpha(v_i)=v_{i+1}$
is an automorphism of $B(m;R,S,T)$, having two vertex-orbits of the same size. 

The circulant graph induced on the set $V_1$ is called the \emph{outer rim} and the circulant graph induced on the set $V_2$ is called the \emph{inner rim}.
We call the vertices from $V_1$ the \emph{outer vertices} and the vertices from $V_2$ the \emph{inner vertices}. There are three types of edges: the edges adjacent to two outer vertices are called \emph{outer edges}, the edges adjacent to two inner vertices are called \emph{inner edges}, and edges connecting an outer vertex to an inner vertex are called  \emph{spokes}. Specifically, the edges $u_iu_{i+a}$, $i \in \ZZ_m$, $a \in R$, are called \emph{outer edges of type} $a$, the edges $v_iv_{i+b}$, $i \in \ZZ_m$, $b \in T$, are called \emph{inner edges of type} $b$ and the edges $u_iv_{i+c}$, $i \in \ZZ_m$, $c \in S$, are called \emph{spokes of type} $c$. 
We will also say that a path is \emph{outer} ( \emph{inner}) if all of its vertices are outer (inner) vertices.

In accordance with our previous discussion, we have $s = |S|$. The order of a graph $B(m;R,S,T)$ is $n = 2m$, the valence is $d$ and $|R| = |T| = d-s$. In the study of bicirculants, other authors use similar notation, see for instance \cite{quarticet}. 

We have already mentioned that generalized rose window graphs are bicirculants. For their description, we need four parameters.
Let $m \ge 3$ be a positive integer and $a,b,c \in \ZZ_m  \setminus\{0\}$  with $a,b \ne m/2$. If we take $R = \{a,-a\}$, $S = \{0,c\}$ and $T = \{ b, -b\}$, the graph $B(m;R,S,T)$ is a generalized rose window graph, which we will denote by $R(m;a,b,c)$.  If $a = 1$, an ordinary rose window graph is obtained.

\begin{figure}[h]
    \centering
\includegraphics[scale=0.5,angle=0]{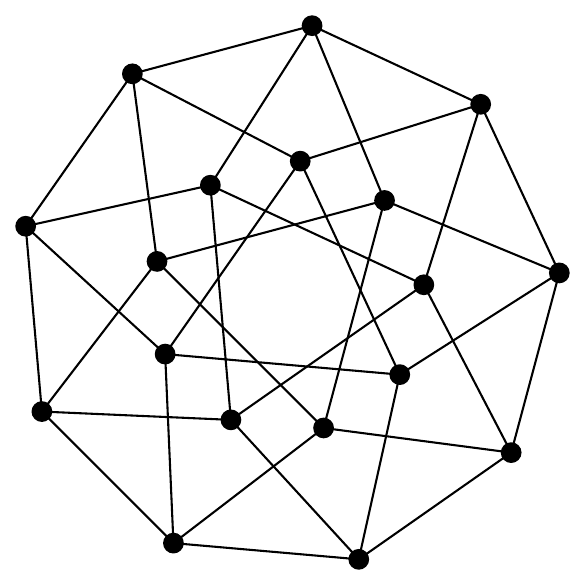}    
\hspace*{2cm}
\includegraphics[scale=0.5,angle=0]{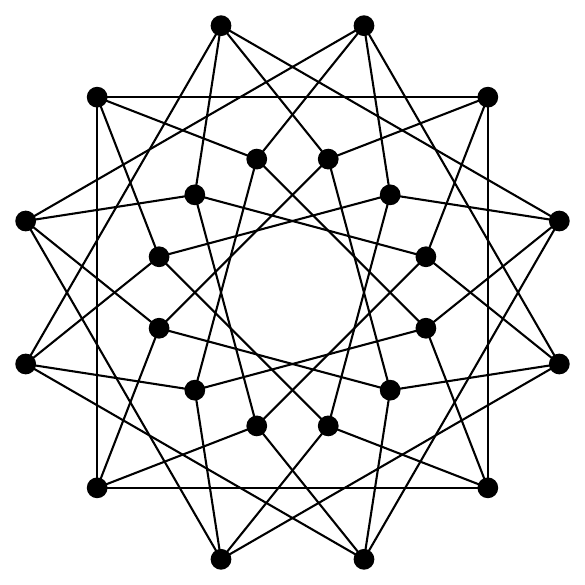}
   \caption{The   generalized rose window graphs $R(9;1,3,2)$ and    $R(12;3,4,2)$.}
    \label{fig:RW9123}
\end{figure}

Figure \ref{fig:RW9123} shows  two generalized rose window graphs. The generalized rose window graph $R(12;3,4,2)$ that is presented on the right hand side of the figure is not isomorphic to any rose window graph.

An $I$-graph $I(m;a,b)$ is a bicirculant $B(m;R,S,T)$ with $m \ge 3$, $R=\{a,-a\}$, $T=\{b,-b\}$  and $R=\{0\}$, where $a,b \in \ZZ_m \setminus \{0,m/2\}$. 
Generalized Petersen graphs are a subfamily of $I$-graphs; an $I$-graph is isomorphic to a generalized Petersen graph if and only if $\gcd(m,a)=1$ or $\gcd(m,b)=1$. We denote the generalized Petersen graph $I(m;1,k)$ by $G(m,k)$.

As we can see, we keep in the description of a specific family of bicirculants for each pair of parameters $x,-x$ only one parameter. Also, we leave out parameters having constant values, such as 0 or 1.

Some properties of bicirculant graphs can be deduced from the general theory of covering graphs \cite{GT}. We will use the following notation. Let $A$ be a set and let $i$ be an integer. We define $A-i=\{a-i \ | \ a \in A \}$ and $A/i= \{a/i \ | \ a \in A \}$.

\begin{proposition}  \label{prop:connected}
A bicirculant   $B(m;R,S,T)$  is connected if and only if $\gcd(m,R,S,T)=1$. In particular, the generalized rose window graph $R(m;a,b,c)$ is connected if and only if $\gcd(m,a,b,c)=1$.
\end{proposition}

In the case where a bicirculant graph is disconnected, it is composed of isomorphic connected components.

\begin{proposition}\label{pro:d_connected_components} 
Let $G = B(m;R,S,T)$. Suppose $\delta = \gcd(m,R,S,T)>1$. 
Then $G$ is a disjoint union of $\delta$   isomorphic graphs $G_0, \dots G_{\delta-1}$ such that $u_i \in G_i$, $i = 0,\dots,\delta-1$. Moreover, each $G_i$  is connected and isomorphic to the graph $B(m/\delta;R/\delta,S/\delta,T/\delta)$. 
\end{proposition}

In many cases there exist isomorphic bicirculants with different parametric descriptions. 
Two special cases are presented below.

\begin{proposition}\label{pro:subtract_isomorphism} 
Graph $B(m;R,S,T)$ is isomorphic to the graph $B(m;R,S-c,T)$ for every $c \in S$.
\end{proposition}

\begin{proposition}\label{pro:arithmetic_isomorphism} 
Let $G = B(m;R,S,T)$, let $r \in \ZZ_m$ be such that $ \gcd(m,r) = 1$ and let $G' = B(m;rR,rS,rT)$. Then the graph $G'$ is isomorphic to the graph $G$.
\end{proposition}

For example, this property of bicirculants was applied to $I$-graphs in the proof that all generalized Petersen graphs are unit-distance graphs \cite{ZiHoPi2012}.
This property also implies that a generalized rose window graph $R(m;a,b,c)$ is isomorphic to  a rose window graph if $\gcd(m,a)=1$ or $\gcd(m,b)=1$.

\section{Hamilton cycles in I-graphs}
\label{sec:Igraphs}

In this section, we consider Hamilton cycles in $I$-graphs as they will play an essential role in the construction of  Hamilton cycles in  rose window graphs. Recall that by removing a set of spokes of the same type  from a rose window graph, we obtain an $I$-graph. We classify Hamilton cycles of $I$-graphs into three types. For each of these types we define a different construction in Section \ref{sec_rosewin}, which shows how to combine  Hamilton cycles in connected components of a rose window graph to a Hamilton cycle in the whole graph.

Hamilton cycles and paths in $I$-graphs are guaranteed by the following results. In \cite{Al11983} Brian Alspach showed that every generalized Petersen graph is hamiltonian, except for the family $G(m,2)$ with $ m \equiv 5 \pmod 6$. However, Alspach and Liu showed that all these exceptional graphs have very many Hamilton paths \cite[Theorem 4.2]{AlspachLiu}.

\begin{theorem} [\cite{AlspachLiu}] \label{thm:GPhamilton}
Every pair of non-adjacent vertices in $G(m,2)$ with $ m \equiv 5 \pmod 6$ is connected by a Hamilton path.
\end{theorem}

Later, it was shown by Bonvicini and Pisanski \cite{BoPi2017} that the non-hamiltonian generalized Petersen graphs are the only non-hamiltonian connected $I$-graphs. Consequently, we have the following theorem.

\begin{theorem}[\cite{BoPi2017}] \label{th_cubicbicirculants}
Every connected $I$-graph, except for the generalized Petersen graphs $G(m,2)$ with $ m \equiv 5 \pmod 6$, is hamiltonian.
\end{theorem}

Clearly, a Hamilton cycle in an $I$-graph $I(m;a,b)$ alternates between paths in the outer rim and paths in the inner rim, which are connected by the spokes. The paths in each rim cover all the vertices of the rim and there are no paths of length zero. This follows from the fact that every vertex of an $I$-graph is adjacent to exactly one spoke. If all the rim paths of a Hamilton cycle $C$ have length one, then $C$ contains all the spokes of the $I$-graph and the spokes alternate with the inner/outer edges; in this case we say that the Hamilton cycle $C$ is \emph{alternating}. Otherwise, it is called \emph{non-alternating}. 

Non-alternating Hamilton cycles are further divided into two types, the $4$-hooked and the $2$-hooked Hamilton cycles. See Section \ref{sec_rosewin} for an explanation of these terms. 
If there exists a labeling of the vertices of the graph $I(m;a,b)$ such that a non-alternating Hamilton cycle $C$ contains the edges $u_0\,u_a$, $u_b\,u_{a+b}$, $v_0\,v_b$,  $v_a\,v_{a+b}$, then $C$ is called \emph{$4$-hooked}. If there exists a labeling of the vertices of the graph $I(m;a,b)$  such that a non-alternating Hamilton cycle $C$ provides a Hamilton path connecting the vertices $v_0$ and $v_a$, or a Hamilton path connecting the vertices $u_0$ and $u_b$, then $C$ is called \emph{$2$-hooked}. By saying that the Hamilton cycle provides a certain Hamilton path, we mean that starting from the cycle, one can produce the path by replacing one or more of its edges with edges not in that cycle. 

Observe that by symmetry,  by adding the same number to the subscripts of the vertices of any given Hamilton cycle, we again obtain a (usually different)   Hamilton cycle. This fact will play a key role in Section \ref{sec_rosewin}.

%
Lemma \ref{lemma_labeling} gives the classification of Hamilton cycles in an $I$-graph $I(m;a,b)$ when $a \ne \pm b$. We deal with the case when $a=b$ or $a=-b$ separately.

\begin{lemma} \label{lem_classification_ham_cycle_Igraph}
Let $G=I(m;a,b)$ be a connected $I$-graph, with $a = b$ or $a=-b$. Then $G$ contains a $2$-hooked Hamilton cycle.
\end{lemma}
\begin{proof}
Observe that $\gcd(m,a)=\gcd(m,b)=1$ since $G$ is connected. Therefore, the sequence $v_0,u_0,u_a,u_{2a},\dots,u_{(m-1)a},$ $v_{(m-1)a},\dots,v_{2a},$ $v_a,v_0$ defines a non-alternating Hamilton cycle, say $C$. By removing the edge $v_0v_a$ from $C$, we obtain a Hamilton path from $v_0$ to $v_a$. That  means that the graph $G$ contains a $2$-hooked Hamilton cycle.
\end{proof}

\begin{lemma}\label{lemma_labeling}
Let $G=I(m;a,b)$ be a connected $I$-graph, with $a \ne \pm b$. Then every Hamilton cycle of $G$ is alternating or $4$-hooked or $2$-hooked.
\end{lemma}

\begin{proof}
Let $C$ be a Hamilton cycle in the graph $G$. Then it is either alternating or non-alternating. We  assume that the cycle $C$ is non-alternating and we will show that it is either $4$-hooked or $2$-hooked.

To this end, we define a special type of non-alternating Hamilton cycles, that we  call almost alternating, and then  we treat separately the cases in which the Hamilton cycle $C$ is almost alternating, and when it is not.  
The Hamilton cycle $C$ is said to be \emph{almost alternating},
if all of the outer and inner subpaths of the   Hamilton cycle $C$ consist of at most two edges and there exists at least one outer or inner subpath with exactly two edges (so the cycle is not alternating). 
\medskip

\noindent
\textbf{Case 1: the Hamilton cycle $C$ is almost alternating.}
Assume that the Hamilton cycle $C$ is almost alternating.
Then there exists at least one outer subpath consisting of two edges. We can label the three vertices of the  subpath with $(u_{-a},  u_0,  u_{a})$ and find the subpath $(v_{-b}, v_0, v_b)$ in $C$ accordingly.  We thus have the edges in $u_0\,u_a$,  $v_0\,v_b$ in $C$. If the edges $u_b\,u_{a+b}$,  $v_a\,v_{a+b}$ are also in $C$, then 
the cycle $C$ is $4$-hooked and the assertion follows. Otherwise we consider two cases: (a) none of the edges $u_b\,u_{a+b}$,  $v_a\,v_{a+b}$ is in $C$ and (b) exactly one of the edges $u_b\,u_{a+b}$,  $v_a\,v_{a+b}$ is in $C$.\\[-3mm]

\textbf{Case (a)} Assume that none of the edges $u_b\,u_{a+b}$,  $v_a\,v_{a+b}$ is in $C$. Then $C$ has the edges $u_b\,u_{b-a}$ and $v_a\,v_{a-b}$ that, combined with the fact that the inner and outer subpaths of $C$ consist of at most three vertices, imply the existence of the subpaths $(u_{b-a}, u_b, v_b, v_0, v_{-b})$ and $(v_{a-b}, v_{a}, u_a, u_0, u_{-a}, v_{-a})$ in $C$. 

If $v_{-a}\,v_{b-a}$ is an edge of $C$, then by adding $a$ modulo $m$ to the subscripts of the vertices of $G$ we see that the cycle $C$ is $4$-hooked, since the edges $v_0\,v_{b}$, $u_b\,u_{b-a}$, $u_0\,u_{-a}$, $v_{-a}\,v_{b-a}$ are turned into the edges $v_a\,v_{a+b}$, $u_{b}\,u_{a+b}$, $u_0\,u_{a}$, $v_{0}\,v_{b}$.

If the edge $v_{-a}\,v_{b-a}$ is not in $C$, then we find the subpaths $(v_{b-a}, u_{b-a}, u_b, v_b, v_0, v_{-b})$ and $(v_{a-b}, v_{a}, u_a, u_0, u_{-a}, v_{-a}, v_{-a-b})$ in $C$. Two cases can occur: the edge $u_{-b}\,u_{-b-a}$ is in $C$ or not.

If $u_{-b}\,u_{-b-a}$ is in $C$, then by adding $(a+b)$ modulo $m$ to the subscripts of the vertices of $G$ we see that the cycle $C$ is $4$-hooked, since the edges $u_{-b}\,u_{-b-a}$, $v_0\,v_{-b}$,  $u_0\,u_{-a}$, $v_{-a}\,v_{-a-b}$ are turned into the edges $u_{0}\,u_{a}$, $v_a\,v_{a+b}$,  $u_b\,u_{a+b}$, $v_{0}\,v_{b}$.

If the edge $u_{-b}\,u_{-b-a}$ is not in $C$, then the edge $u_{-b}\,u_{a-b}$ is in $C$ and by adding $b$ modulo $m$ to the subscripts of the vertices of $G$, the edges $u_{-b}\,u_{a-b}$, $v_0\,v_{-b}$,  $u_0\,u_{a}$, $v_{a}\,v_{a-b}$ are turned into the edges $u_{0}\,u_{a}$, $v_0\,v_{b}$,  $u_b\,u_{a+b}$, $v_{a}\,v_{a+b}$, therefore $C$ is again a $4$-hooked Hamilton cycle.

\textbf{Case (b)} Now assume that exactly one of the edges $u_b\,u_{a+b}$,  $v_a\,v_{a+b}$ is in $C$.
For the case where $u_b\,u_{a+b}$ is in $C$ but $v_a\,v_{a+b}$ is not in $C$, we can repeat the same arguments as above when the edge $u_{-b}\,u_{a-b}$ is in $C$, and also when both edges $u_{-b}\,u_{-a-b}$, $v_{-a}\,v_{-a-b}$ are in $C$. In the missing case, that is, when 
$u_{-b}\,u_{-a-b}$ and $v_{-a}\,v_{b-a}$ are edges of $C$, we can find a Hamilton path from $u_0$ to $u_b$, or from $v_0$ to $v_a$ in $G$, so the cycle $C$ is $2$-hooked.
More specifically, we consider the vertices in clockwise order, and we can always assume that $v_b$, $v_0$ occur in that order. The vertices occur in $C$ in one of the following orders: $v_{a+b}$, $u_{a+b}$, $u_b$, $v_b$, $v_0$, $v_{a-b}$, $v_a$, $u_a$, $u_0$, $u_{-a}$, 
or $v_{a+b}$, $u_{a+b}$, $u_b$, $v_b$, $v_0$, $u_{b-a}$, $v_{b-a}$, $v_{-a}$, $u_{-a}$, $u_0$, $u_{a}$.
In the first case we remove the edges $u_a\,v_a$, $u_{a+b}\,v_{a+b}$, add the chord $v_a\,v_{a+b}$, and find a Hamilton path from $u_a$ to $u_{a+b}$ that yields a Hamilton path from $u_0$ to $u_b$ if we add $(-a)$ modulo $m$ to the subscripts of the vertices of $G$; see the graph on the left-hand side of Figure \ref{fig_labeling_almost_alternating_cycles}.
In the second case we remove the edges $u_{b}\,v_{b}$, $u_{b-a}\,v_{b-a}$, add the chord $u_b\,u_{b-a}$, and find a Hamilton path from $v_{b-a}$ to $v_{b}$ that yields a Hamilton path from $v_0$ to $v_a$ if we add $(a-b)$ modulo $m$ to the subscripts of the vertices of $G$; see the graph on the right-hand side of Figure \ref{fig_labeling_almost_alternating_cycles}.

\begin{figure}[htb]
\begin{center}
\includegraphics[width=12cm]{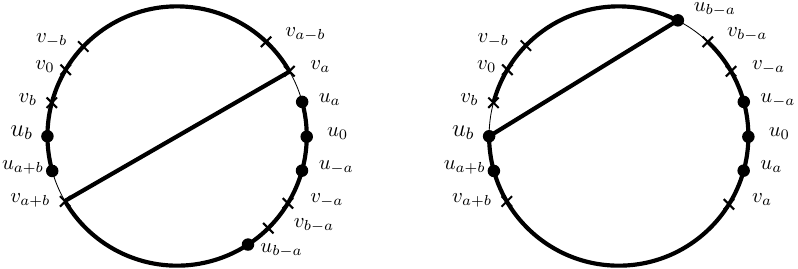}
\end{center}
\caption{
The almost alternating Hamilton cycle $C$ described in Lemma \ref{lemma_labeling} when the edge $u_b\,u_{a+b}$ is in $C$, but the edge $v_a\,v_{a+b}$ is not. 
The bold lines define a Hamilton path from $u_a$ to $u_{a+b}$ in the cycle on the left-hand side of the figure and a Hamilton path from $v_b$ to $v_{b-a}$ in the cycle on the right-hand side.}
\label{fig_labeling_almost_alternating_cycles}
\end{figure}

For the case where $v_a\,v_{a+b}$ is in $C$ but  $u_b\,u_{a+b}$ is not in $C$ we can repeat the same arguments as for the case where $u_b\,u_{a+b}$ is in $C$ but  $v_a\,v_{a+b}$ is not in $C$  by symmetry.
The validity of the lemma is thus proved for the almost alternating cycles of an $I$-graph.
\medskip

\noindent
\textbf{Case 2: the Hamilton cycle $C$ is 
not almost alternating. }
Assume now that the Hamilton cycle $C$ is not almost alternating.
Then we find at least one outer or inner subpath in $C$ consisting of at least three edges; by symmetry we may assume that such a path is an outer subpath and we can label the four vertices of the subpath with $(u_{-a}, u_0, u_a, u_{2a})$. We find the inner subpaths $(v_{-b}, v_0, v_{b})$,  $(v_{a-b}, v_a, v_{a+b})$ in $C$ accordingly. 
We therefore find the edges $v_0\,v_b$, $u_0\,u_a$, $v_a\,v_{a+b}$ in $C$. 
If the edge $u_b\,u_{a+b}$ is also in $C$, then the cycle $C$ is $4$-hooked and the assertion follows.

We now assume that the edge $u_b\,u_{a+b}$ is not in $C$. Then the Hamilton cycle $C$ is $2$-hooked -- and hence the assertion follows -- if the vertices occur  in $C$ in some prescribed orders. More specifically, in the following we consider the vertices of $C$ in clockwise order; we can always assume that the vertices $u_a$, $u_0$ occur in that order; we also set $\{v_x, v_{x’}\}=\{v_{a+b}, v_{a-b}\}$,  $\{v_y, v_{y’}\}=\{v_{b}, v_{-b}\}$.

If the vertices occur in order $u_a$, $u_0$, $v_x$, $v_a$, $v_{x’}$, $v_y$, $v_0$, $v_{y’}$, we can find a Hamilton path from $v_0$ to $v_a$ in $G$. In fact, at least one of the equalities $x-y =a$ or $x-y'=a$ holds. If $x-y=a$ (respectively, $x-y’=a$), then we remove the edges $u_0\,u_a$, $v_a\,v_x$, $v_0\,v_{y}$ (respectively, $u_0\,u_a$, $v_a\,v_x$, $v_0\,v_{y'}$) in $C$, add the chords $u_0\,v_0$, $u_a\,v_a$, and find a Hamilton path from $v_x$ to $v_{y}$ (respectively, from $v_x$ to $v_y'$)  that yields a Hamilton path from $v_0$ to $v_a$, since $x-y=a$ (respectively, $x-y’=a$); see the first two graphs of Figure \ref{fig_labeling_almost_alternating_cycles_2}.
We can find a Hamilton path from $v_0$ to $v_a$ even in the case where the vertices occur in the order $u_a$, $u_0$, $v_y$, $v_0$, $v_{y’}$, $v_x$, $v_a$, $v_{x’}$, with $x’-y=a$. In fact, in this case we remove the edges $u_0\,u_{a}$, $v_a\,v_{x’}$, $v_0\,v_{y}$ in $C$, add the chords $u_0\,v_0$, $u_a\,v_a$, and find a Hamilton path from $v_y$ to $v_{x’}$  that yields a Hamilton path from $v_0$ to $v_a$, since $x’-y=a$;  see the third graph of Figure \ref{fig_labeling_almost_alternating_cycles_2}.

\begin{figure}[htb]
\begin{center}
\includegraphics[width=16cm]{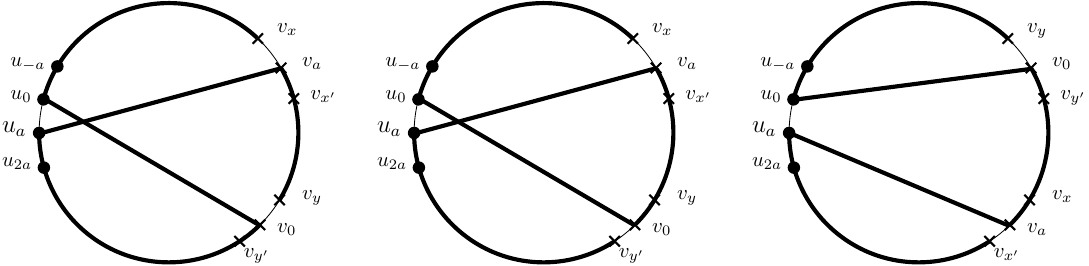}
\end{center}
\caption{
The Hamilton cycle $C$ described in Lemma \ref{lemma_labeling}:
the cycle is not almost alternating and the vertices occur in the order $u_a$, $u_0$, $v_x$, $v_a$, $v_{x’}$, $v_y$, $v_0$, $v_{y’}$  with $\{v_x, v_{x’}\}=\{v_{a+b}, v_{a-b}\}$,  $\{v_y, v_{y’}\}=\{v_{b}, v_{-b}\}$  in the first two cycles and in the order $u_a$, $u_0$, $v_y$, $v_0$, $v_{y’}$, $v_x$, $v_a$, $v_{x’}$ with $x'-y=a$ in the third cycle. The bold lines in the first two cycles define a Hamilton path from $v_x$ to $v_y$, and a Hamilton path from $v_x$ to $v_{y'}$. The bold lines in the third cycle define a Hamilton path from $v_y$ to $v_{x'}$.}\label{fig_labeling_almost_alternating_cycles_2}
\end{figure}

\begin{figure}[htb]
\begin{center}
\includegraphics[width=12cm]{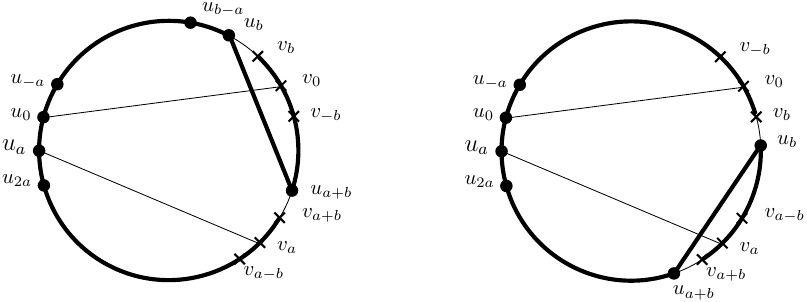}
\end{center}
\caption{
The Hamilton cycle $C$ described in Lemma \ref{lemma_labeling} when
the cycle is not almost alternating, does not contain the edge  $u_b\,u_{a+b}$, and the vertices occur in the order $u_a$, $u_0$, $v_y$, $v_0$, $v_{y’}$, $v_x$, $v_a$, $v_{x’}$ with $x'-y\neq a$. Furthermore, $v_y=v_b$ and $v_{x'}=v_{a-b}$ in the cycle on the left-hand side and $v_y=v_{-b}$ and $v_{x'}=v_{a+b}$ in the cycle
on the right-hand side. The bold lines define a Hamilton path from $v_b$ to $v_{a+b}$.}\label{fig_labeling_almost_alternating_cycles_3}
\end{figure}

It remains to consider the case in which the edge $u_b\,u_{a+b}$ is not in $C$, and the vertices occur in the order   $u_a$, $u_0$, $v_y$, $v_0$, $v_{y’}$, $v_x$, $v_a$, $v_{x’}$, with $x’-y\neq a$, i.e., $(v_y, v_{x'})=(v_b, v_{a-b})$ or $(v_y, v_{x'})=(v_{-b}, v_{a+b})$. The nonexistence of the edge $u_b\,u_{b+a}$ in $C$ implies the existence of the subpaths $(u_{b-a}, u_{b}, v_b, v_0, v_{-b})$ and $(u_{b+a}, v_{b+a}, v_{a}, v_{a-b})$ in $C$. In this setting, we remove the edges $u_b\,v_b$, $u_{a+b}\,v_{a+b}$ from $C$, and add the edge $u_b\,u_{a+b}$; we find a Hamilton path from $v_b$ to $v_{a+b}$ that provides a Hamilton path from $v_0$ to $v_a$ if we add $-b$ modulo $m$ to the subscripts of the vertices of $G$; see Figure \ref{fig_labeling_almost_alternating_cycles_3}. 
Therefore, the cycle $C$ is $2$-hooked and the assertion follows. 
\end{proof}

When constructing a Hamilton cycle  in a rose window graph using a $4$-hooked Hamilton cycle in its subgraph that is an $I$-graph, certain orderings of vertices are difficult to deal with.
A $4$-hooked Hamilton cycle $C$  is called \emph{elusive (of type $1$ or of type $2$)}, if the vertices $u_0,u_a, u_b,u_{a+b},v_0,v_b,v_a,v_{a+b}$ appear on $C$ in one of the orders \eqref{order1} or \eqref{order2}, starting with vertices $u_0,u_a$:
\begin{align}  
 & u_0,u_a,u_b,u_{a+b},v_{a+b},v_a,v_b,v_0,   
 \label{order1} \\  
 & u_0,u_a,v_a,v_{a+b},v_0,v_b,u_b,u_{a+b}.  \label{order2}
\end{align}
Otherwise, a $4$-hooked Hamilton cycle is called \emph{standard}. 
Both types of elusive Hamilton cycles are equivalent in a certain way.

\begin{remark} \label{remark1_2}
Let $C$ be an elusive Hamilton cycle of type 2 in an $I$-graph $I(m;a,b)$. We relabel the vertices by adding $-a$ modulo $m$ to their indices. Then the sequence of vertices $u_0,u_a,v_a,v_{a+b},v_0,v_b,u_b,u_{a+b}$ is mapped to the sequence 
$u_{-a},u_0,v_0,v_{b},v_{-a}, v_{-a+b},$ $ u_{-a+b},u_{b}$. By reversing the cycle $C$ we see that it contains the sequence $u_0,u_{-a},u_b,   u_{-a+b},$ $v_{-a+b},v_{-a},$ $v_b,v_0$, so $C$ is an elusive Hamilton cycle of type 1 for the graph $I(m;-a,b)$, which is the same graph as $I(m;a,b)$. 
\end{remark}

\begin{remark} \label{remark_lemma8}
One can observe that  an elusive Hamilton cycle $C$ that may appear  in the proof of Lemma \ref{lemma_labeling}  is either  of type 1 with  vertices $u_0,v_0$ not adjacent in $C$ or it is of type 2 with the property that  it  contains the subpaths $(v_{a+b}, v_a, u_a, u_{0}, u_{-a},v_{-a})$ and $(u_{a+b},u_b,v_b,v_0,v_{-b},u_{-b})$ 
occurring in this order in $C$; see the first part of the proof regarding the almost alternating Hamilton cycles.
By Remark \ref{remark1_2} we may thus assume that such a cycle is also of type 1 and that it contains the subpaths 
$(v_b, v_0, u_0, u_a, u_{2a},v_{2a})$ 
and $(u_b, u_{a+b},$ $ v_{a+b},$ $v_a,$ $v_{a-b}, u_{a-b})$
occurring in this order in $C$.
\end{remark}

The next lemma shows that it is almost always possible to replace elusive Hamilton cycles of type 1 in $I$-graphs with standard Hamilton cycles or certain Hamilton paths. The proof of the lemma   is rather long and is given in Appendix \ref{appendix}.

\begin{lemma}  \label{lem_4hookedmain}
Let $a \ne \pm b$ and let an $I$-graph $G=I(m;a,b)$ contain an elusive Hamilton cycle $C$ of type 1. Then 
\begin{itemize}
    \item the graph $G$ contains a standard $4$-hooked Hamilton cycle
    or a $2$-hooked Hamilton cycle or
    \item  $b   \equiv -2a\pmod m$ or $a   \equiv -2b\pmod m$
    and the cycle $C$  contains the subpaths $(v_b, v_0, u_0, u_a)$ 
    and $(u_b, u_{a+b},$ $ v_{a+b}, v_a, v_{a-b}, u_{a-b}) $
    occurring in this order in $C$.
\end{itemize}
\end{lemma}

\section{Hamilton cycles in generalized rose window graphs}
\label{sec_rosewin}

In this section, we show how to construct a Hamilton cycle in any given generalized rose window graph. We will use the following notation: $G=R(m;a,b,c)$ will denote a connected generalized rose window graph, so $\gcd(m,a,b,c)=1$. By $H$ we denote the graph obtained from $G$ by removing the spokes of type $c$. 
Note that the graph $H$ is composed of $\gcd(m, a, b)$ isomorphic connected $I$-graphs.
If the graph $H$ is connected, then Theorem \ref{th_cubicbicirculants} implies $G$ is hamiltonian in case it is not isomorphic to a generalized Petersen graph $G(n, 2)$, $n\equiv 5\pmod 6$; this case has to be considered separately. 

We now consider the case where $H$ is not connected. Set $\lambda=\gcd(m, a,  b)-1$ and denote by $H_i$, $0\leq i\leq\lambda$,  the connected components of $H$; $H_0$ will be the component containing the vertex $u_0$.

The connected component $H_i$ with $i>0$ can be described as the $i$-th isomorphic copy of $H_0$, that is, we leave invariant the adjacencies in $H_0$ and label the vertices of $H_i$ by adding $i\,c$ modulo $m$ to the subscripts of the vertices in $H_0$.
We will use the notation $u^i_x$, $v^i_x$ to denote the outer and inner vertices of $H_i$ corresponding to the outer and inner vertices $u_x$, $v_x$, respectively, in $H_0$.That is, $u^i_x=u_{x+i\,c}$ and $v^i_x=v_{x+i\,c}$. The outer vertices $u^i_x$ in $H_i$ are adjacent to the inner vertices $v^{i+1}_{x}$ in $H_{i+1}$, since $G$ is connected and $H$ is obtained from $G$ by removing the spokes of type $c$. 
Sometimes, for our convenience, we will also use the notation $u^0_x$, $v^0_x$ for the vertices in $H_0$.

\smallskip

Given a generalized rose window graph $G$ whose subgraph $H$ has at least two connected components, we will construct a Hamilton cycle in $G$ by appropriately joining the Hamilton cycles, or paths, in the components $H_i$.  The construction depends on the classification defined in Lemma \ref{lemma_labeling}. More specifically, for an alternating Hamilton cycle we will define the \emph{alternating construction} (see Proposition \ref{pro_alternating_construction}); 
we will define the $4$- and the $2$-\emph{hooked construction} for the $4$- and the $2$-hooked Hamilton cycles, respectively. 
The terminology follows from the fact that, in the assembly of the cycles in the components of $H$, the cycle corresponding to $C$ in $H_i$ with $0\leq i\leq\lambda-1$, will be connected to the cycle corresponding to $C$ in $H_{i+1}$ by $4$ or $2$ spokes, respectively.

\begin{figure}
\begin{center}
\includegraphics[width=12cm]{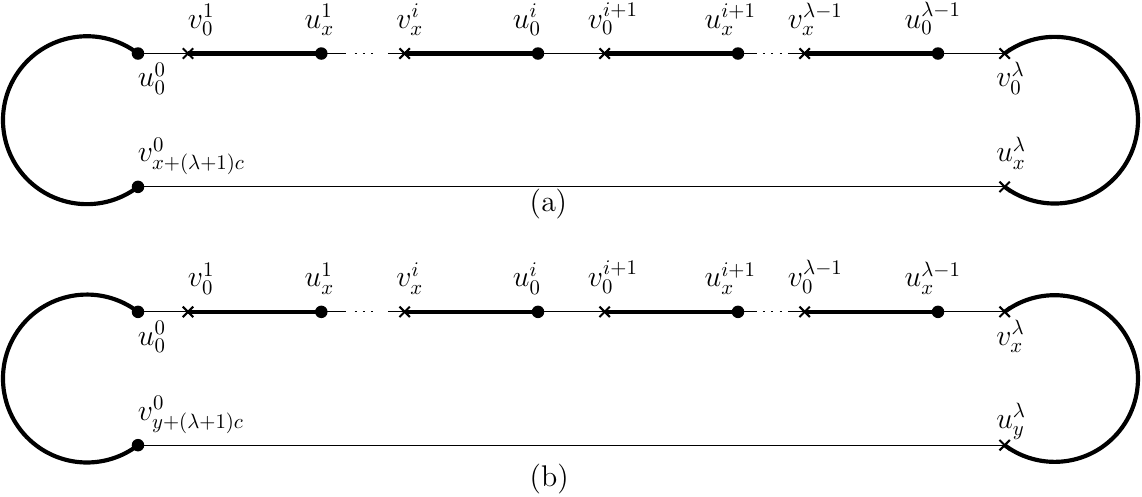}
\end{center}
\caption{The construction of a Hamilton cycle in a generalized rose window graph $G$ when the connected components of $H$ are isomorphic to a generalized Petersen graph $G(n, 2)$, $n\equiv 5\pmod 6$; see the proof of Theorem \ref{th_rosewin_hamiltonian}.  We apply case (a) for odd values of $\lambda=\gcd(m, a, b)-1$, and case (b) for even values of $\lambda$. 
The bold lines represent the paths $v^i_0\,P\,u^i_x$,  $v^i_x\,P\,u^i_0$ with $1\leq i\leq\lambda-1$, the Hamilton paths $u^0_0\,P\,v^0_{x+(\lambda+1)c}$
and $v^{\lambda}_0\,P\,u^{\lambda}_{x}$ in case (a) and  the Hamilton paths
 $u^0_0\,P\,v^0_{y+(\lambda+1)c}$ and  $v^{\lambda}_x\,P\,u^{\lambda}_{y}$ in case (b).  }\label{fig_generalized_petersen}
\end{figure}

The alternating and the hooked constructions can be applied when $H_0$ is not isomorphic to the generalized Petersen graph $G(n, 2)$ with $n\equiv 5\pmod 6$. In the latter case, the graph does not have a Hamilton cycle and we will apply the construction described in the proof of Theorem \ref{th_rosewin_hamiltonian}, and summarized in Figure \ref{fig_generalized_petersen}, which could be called the $1$-hooked construction in analogy to the previous ones. 

We now give the alternating and the hooked constructions; the $2$-hooked construction will be also used in Proposition \ref{pro_4hooked}, which defines the $4$-hooked construction.  In what follows, we will use the notation $x\, P\, y$ to denote a path $P$ from the vertex $x$ to the vertex $y$.

\begin{proposition}\label{pro_alternating_construction}{\bf The alternating construction.}
Let $G=R(m;a,b,c)$ be a connected generalized rose window graph.
Let $H$ be the graph obtained from $G$ by removing the spokes of type $c$,  and let $H_0$ be the connected component of $H$ containing the vertex $u_0$.
Assume $\lambda=\gcd(m, a, b)-1> 0$.
If the graph $H_0$ has an alternating Hamilton cycle, then the graph $G$ is hamiltonian.
\end{proposition}

\begin{proof}  Let  $m_0=m/\gcd(m,a,b)$. Assume that the graph $H_0$ has an alternating Hamilton cycle; denote it by $C$. Note that the existence of an alternating Hamilton cycle   in $H_0$ implies that $m_0$ is even. 
We denote the outer and inner vertices of $H_0$ with $u_{x_j}$, $v_{x_j}$, respectively, for $0\leq j\leq m_0-1$, so that the vertices $u_{x_j}$,  $v_{x_j}$ are consecutive in $C$, as well as $u_{x_j}$, $u_{x_{j+1}}$ and $v_{x_{j-1}}$, $v_{x_j}$,  with $1\leq j\leq m_0-1$, $j$ odd. The indices $x_j$ are integers modulo $m$.  
In $H_i$ with $i>0$,  the vertices corresponding to $u_{x_j}$, $v_{x_j}$ of $H_0$ will be denoted with $u^i_{x_j}$, $v^i_{x_j}$; the vertices in $H_0$
will be also denoted with $u^0_{x_j}$, $v^0_{x_j}$.

To construct a Hamilton cycle in the graph $G$, we keep just the spokes $v^i_{x_j}\,u^i_{x_j}$ (spokes of type 0) in each of the graphs $H_i$ for  $1\leq i\leq\lambda-1$.
For every $0\leq j\leq m_0-1$,  we connect  the edge $v^i_{x_j}\,u^i_{x_j}$ in $H_i$ to the edge $v^{i+1}_{x_j}\,u^{i+1}_{x_j}$ in $H_{i+1}$ by adding the spoke $u^i_{x_j}\,v^{i+1}_{x_j}$ of type $c$, for  $1\leq i\leq\lambda-2$. 
For every $0\leq j\leq m_0-1$,  we  get a path from $v^1_{x_j}$ to $u^{\lambda-1}_{x_j}$, to which we add the edges $u^0_{x_j}\,v^1_{x_j}$ and $u^{\lambda-1}_{x_{j}}\,v^{\lambda}_{x_j}$ in order to obtain a path $u^0_{x_j}\,P\,v^{\lambda}_{x_j}$ from the vertex $u^0_{x_j}$ in $H_0$ to the vertex $v^{\lambda}_{x_j}$ in $H_{\lambda}$.  The union of the paths $u^0_{x_j}\,P\,v^{\lambda}_{x_j}$
is a disconnected graph that covers all the vertices in $G$, with the exception for the inner vertices in $H_0$ and the outer vertices in $H_{\lambda}$.
Since $m_0$ is even,  we can join the paths $u^0_{x_j}\,P\,v^{\lambda}_{x_j}$ by adding the paths $(u^0_{x_j}, v^0_{x_j}, v^0_{x_{j+1}}, u^0_{x_{j+1}})$ for $0\leq j\leq m_0-1$, $j$ even,
and the paths $(v^{\lambda}_{x_j}, u^{\lambda}_{x_j}, u^{\lambda}_{x_{j+1}}, v^{\lambda}_{x_{j+1}})$ for $0\leq j\leq m_0-1$, $j$ odd.  
We thus obtain a Hamilton cycle in $G$.  We summarize the construction with the diagram in Figure \ref{fig_alternating_construction}.
\begin{figure}[htb]
\begin{center}
\includegraphics[width=12cm]{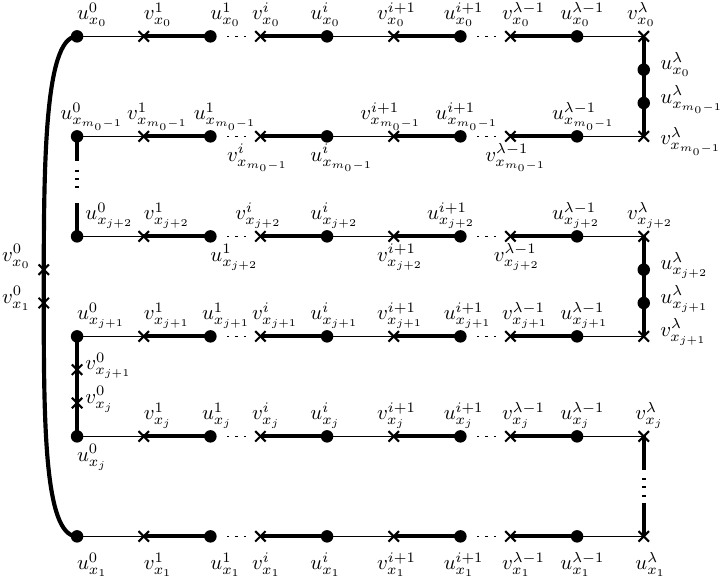}
\end{center}
\caption{The alternating construction for the generalized rose window graphs described in Proposition \ref{pro_alternating_construction}. The bold lines represent
the edges $v^i_{x_j}\,u^i_{x_j}$ in $H_i$,  $1\leq i\leq\lambda-1$, the paths $(u^0_{x_j}, v^0_{x_j}, v^0_{x_{j+1}}, u^0_{x_{j+1}})$ for $0\leq j\leq m_0-1$, $j$ even,
and the paths $(v^{\lambda}_{x_j}, u^{\lambda}_{x_j}, u^{\lambda}_{x_{j+1}}, v^{\lambda}_{x_{j+1}})$ for $0\leq j\leq m_0-1$, $j$ odd.}
\label{fig_alternating_construction}
\end{figure} 
\end{proof}

\begin{proposition}\label{pro_2hooked}{\bf The $2$-hooked construction.}
Let $G=R(m;a,b,c)$ be a connected generalized rose window graph.
Let $H$ be the graph obtained from $G$ by removing the spokes of type $c$,  and let $H_0$ be the connected component of $H$ containing the vertex $u_0$.
Assume $\lambda=\gcd(m, a, b)-1> 0$.
If  the graph $H_0$ has a $2$-hooked Hamilton cycle, then  the graph $G$ is hamiltonian.
\end{proposition}

\begin{proof}
Assume that the graph $H_0$ has a $2$-hooked   Hamilton cycle; denote it by $C$.
We  first assume that the $2$-hooked cycle $C$ provides a Hamilton path connecting the vertices $v_0$ and $v_a$ in $H_0$, say $v_0\,P\,v_a$. Such a path necessarily contains an outer edge, say $u_x\,u_{x+a}$, since it contains the same number of outer and inner vertices.
Without loss of generality, we can assume that $u_x$ precedes $u_{x+a}$ in  $v_0\,P\,v_a$, so that the removal of the edge $u_x\,u_{x+a}$ yields the two subpaths
$v_0\,P\,u_x$ and $v_a\,P\,u_{x+a}$. We can also find a Hamilton path $u_0\,P\,u_a$ from $u_0$ to $u_a$ in $H_0$, since the graph $H_0$ is hamiltonian 
and every Hamilton cycle in $H_0$ has at least one outer edge.

If $\lambda=1$, so the graph $H$ has two components, we connect the two Hamilton paths $u^0_0\,P\,u^0_a$ in $H_0$ and $v^{1}_{0}\,P\,v^{1}_{a}$ in $H_{1}$ by  adding the spokes $u^0_0\,v^1_0$ and $u^0_a\,v^1_{a}$ in order to obtain a Hamilton cycle in $G$.

Now we assume that $\lambda >1$.
In $H_i$, with $1\leq i\leq\lambda-1$,  we consider the subpaths  $v^i_0\,P\,u^i_x$ and $v^i_a\,P\,u^i_{x+a}$  corresponding to the subpaths $v_0\,P\,u_x$ and $v_a\,P\,u_{x+a}$ of $H_0$.  
We turn the subpaths $v^i_0\,P\,u^i_x$ and $v^i_a\,P\,u^i_{x+a}$ into the subpaths $v^i_{(i-1)x}\,P\,u^i_{i\,x}$ and $v^i_{a+(i-1)x}\,P\,  u^i_{a+i\,x}$ by adding $(i-1)x$ modulo $m$ to the subscripts of the vertices in $H_i$.
Notice that by adding $(i-1)x$ modulo $m$ to the subscripts of the vertices in 
$H_i$, we still get vertices of $H_i$.

We now construct a Hamilton cycle in $G$ by connecting the above paths as follows.
For $1\leq i\leq\lambda-2$,  we join the path $v^i_{(i-1)x}\,P\,u^i_{i\,x}$  in $H_i$  to the path $v^{i+1}_{i\,x}\,P\,u^{i+1}_{(i+1)x}$ in $H_{i+1}$  by the spoke $u^i_{i\,x}\,v^{i+1}_{i\,x}$, and
also join the path $v^i_{a+(i-1)x}\,P\,u^i_{a+i\,x}$ in $H_i$ to the path $v^{i+1}_{a+i\,x}\,P\,u^{i+1}_{a+(i+1)x}$ in $H_{i+1}$ by the spoke $u^i_{a+i\,x}\,v^{i+1}_{a+i\,x}$.
We obtain two vertex-disjoint paths -- the former from $v^1_0$ to $u^{\lambda-1}_{(\lambda-1)x}$ and the latter from $v^1_{a}$ to $u^{\lambda-1}_{a+(\lambda-1)x}$ -- 
whose union covers all the vertices in $G-(H_0\cup H_{\lambda})$.  We connect the two paths to the Hamilton paths $u^0_0\,P\,u^0_a$ in $H_0$ and  
$v^{\lambda}_{(\lambda-1)x}\,P\,v^{\lambda}_{a+(\lambda-1)x}$ in $H_{\lambda}$ by adding the spokes $u^0_0\,v^1_0$, $u^0_a\,v^1_{a}$ and 
$u^{\lambda-1}_{(\lambda-1)x}\,v^{\lambda}_{(\lambda-1)x}$, $u^{\lambda-1}_{a+(\lambda-1)x}\,v^{\lambda}_{a+(\lambda-1)x}$. We thus obtain a Hamilton cycle in $G$.
We summarize the construction with the diagram in Figure \ref{fig_rosewin_construction1}.

For the case where the $2$-hooked cycle in $H_0$ provides a Hamilton path connecting $u_0$ and $u_b$, we can repeat the same argument as above (it suffices to
replace the parameter $a$ with the parameter $b$).  The assertion follows.
\begin{figure}[htb]
\begin{center}
\includegraphics[width=14cm]{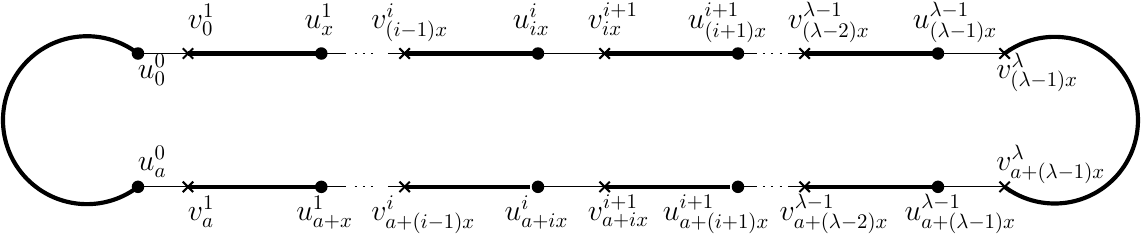}
\end{center}
\caption{The $2$-hooked construction for the generalized rose window graphs described in Proposition \ref{pro_2hooked}. The bold lines represent the paths 
$v^i_{(i-1)x}\,P\,u^i_{i\,x}$,  $v^i_{a+(i-1)x}\,P\,u^i_{a+i\,x}$ in $H_i$ with $1\leq i\leq\lambda-1$, and  the Hamilton paths $u^0_0\,P\,u^0_a$ in $H_0$,  $v^{\lambda}_{(\lambda-1)x}\,P\,v^{\lambda}_{a+(\lambda-1)x}$ in $H_{\lambda}$. The paths are joined by spokes $u^i_{i\,x}\,v^{i+1}_{i\,x}$, $u^i_{a+i\,x}\,v^{i+1}_{a+i\,x}$ for $0\leq i\leq\lambda-1$.\\
}\label{fig_rosewin_construction1}
\end{figure}
\end{proof}

\goodbreak

\begin{remark}\label{rem_2hooked}
     We can apply the $2$-hooked construction described in Proposition \ref{pro_2hooked} even when we have a Hamilton path from $u_0$ to $u_p$ in $H$  (denoted by $u_0\,P\,u_p$), 
    where $u_p$ is an arbitrary vertex of $H$, a Hamilton path from $v_0$ to $v_p$, (denoted by $v_0\,P\,v_p$),  and two paths whose union partitions the vertices of $H$, namely the
    paths $u_0\,P\,v_p$ and $u_p\,P\,v_0$ from $u_0$ to $v_p$ and from $u_p$ to $v_0$, or the paths $u_0\,P\,v_0$ and $u_p\,P\,v_p$
    from $u_0$ to $v_0$ and from $u_p$ to $v_p$, respectively. In fact, in the construction described in the proof of Proposition \ref{pro_2hooked}, we can replace the subpaths $v_0\, P\, u_x$ and $v_a\, P\, u_{x+a}$ in $H_0$ with the subpaths $u_0\,P\,v_p$ and $u_p\,P\,v_0$, or $u_0\,P\,v_0$ and $u_p\,P\,v_p$; consequently, in each $H_i$, with $1\le i\le\lambda-1$, we will consider the corresponding copies; we will take the Hamilton path from $u_0$ to $u_p$ in $H_0$, and the copy of the Hamilton path from $v_0$ to $v_p$ in $H_{\lambda}$. Roughly speaking, in order to obtain a Hamilton cycle in $G$, we will join the outer vertices of $H_i$ to the inner vertices of $H_{i+1}$ having the same subscripts.
    For instance, if we have the paths $u_0\,P\,v_p$ and $u_p\,P\,v_0$ partitioning the vertices in $H_0$, then in $H_i$, with $1\le i\le\lambda-1$, we can consider the corresponding copies $u^i_0\,P\,v^i_p$ and $u^i_p\,P\,v^i_0$ - we recall that $u^i_0=u_{ic}$, $u^i_p=u_{p+ic}$
    and $v^i_0=v_{ic}$, $v^i_{p}=v_{p+ic}$; we find a Hamilton cycle in $G$ by connecting the vertices $u^i_p$ and $u^i_0$, with $1\le i\le\lambda -2$,
    to the vertices  $v^{i+1}_p$ and $v^{i+1}_0$, respectively; we also add the edges $u^0_0v^1_0$, $u^0_pv^1_p$, and $u^{\lambda-1}_0v^{\lambda}_0$,
     $u^{\lambda-1}_pv^{\lambda}_p$.
    \end{remark}

  In Example \ref{example}, we show how to use the $2$-hooked construction, described in Remark \ref{rem_2hooked}. We will need this example in the proof of Proposition \ref{pro_4hooked}.

    \begin{example} \label{example}
    Let $b\equiv -2a\pmod m$, and assume that a Hamilton cycle $C$ of an $I$-graph $H=I(m;a,b)$  contains the subpaths $(v_b, v_0, u_0, u_a)=(v_{-2a}, v_0, u_0, u_a)$ and $(u_b, u_{a+b}, v_{a+b}, v_a, v_{a-b}, u_{a-b})=(u_{-2a}, u_{-a}, v_{-a}, v_{a}, v_{3a}, u_{3a})$ occurring in this order in $C$. Then we find a Hamilton path from $u_0$ to $u_{3a}$, a Hamilton path from $v_0$ to $v_{3a}$, and the paths from $u_0$ to $v_{3a}$ and from $u_{3a}$ to $v_0$, whose union partitions the vertices of $H$. In detail, the existence of the two paths from $u_0$ to $v_{3a}$ and from $u_{3a}$ to $v_0$ is straightforward (remove the edges $u_0v_0$ and $u_{3a}v_{3a}$); the Hamilton path from $v_0$ to $v_{3a}$ can be obtained as follows: remove the edges $v_0v_{-2a}$, $u_0u_a$, $u_{-a}u_{-2a}$, $v_av_{3a}$ from $C$, and add the edges $u_0u_{-a}$, $u_av_a$, $u_{-2a}v_{-2a}$. For the Hamilton path from $u_0$ to $u_{3a}$, we first note that $C$ also contains  the subpath $(u_a,u_{2a},v_{2a},v_{4a})$ and the edge $u_{3a}u_{4a}$.  Then we remove the edges $u_0v_0$,  $v_{2a}v_{4a}$, $u_{3a}u_{4a}$, and add the edges $v_0v_{2a}$, $u_{4a}v_{4a}$.

    Let $a\equiv -2b\pmod m$ and again assume that a Hamilton cycle $C$ of an $I$-graph $H=I(m;a,b)$ contains the subpaths$(v_b, v_0, u_0, u_a)$  and $(u_b, u_{a+b}, v_{a+b}, v_a, v_{a-b}, u_{a-b})$ occuring in this order in $C$.  Notice that in this case the cycle $C$ also contains the subpaths $(u_a,u_0,v_0,v_b,$ $v_{2b}, u_{2b},u_{4b})$ and 
    $(u_{a-b},v_{a-b},v_a,v_{a+b},u_{a+b},u_b,u_{3b},v_{3b},v_{4b})$ in this order and we can repeat arguments similar to the previous ones due to the symmetry between the parameters $a$ and $b$ and find  a Hamilton path from $u_0$ to $u_{3b}$, a Hamilton path from $v_0$ to $v_{3b}$, and the paths from $u_0$ to $v_{3b}$ and from $u_{3b}$ to $v_0$, whose union partitions the vertices of $H$.      
\end{example} 

\begin{proposition}\label{pro_4hooked}{\bf The $4$-hooked construction.}
Let $G=R(m;a,b,c)$ be a connected generalized rose window graph with $a \ne \pm b$.
Let $H$ be the graph obtained from $G$ by removing the spokes of type $c$,  and let $H_0$ be the connected component  of $H$ containing the vertex $u_0$.
Assume $\lambda=\gcd(m, a, b)-1> 0$. 
If $H_0$ has a $4$-hooked Hamilton cycle, then $G$ is hamiltonian.
\end{proposition}

\begin{proof}
Assume that the graph $H_0$ has a $4$-hooked Hamilton cycle; denote it by $C$. Since $C$ is a $4$-hooked cycle, it contains the edges $u_0\,u_a$, $u_b\,u_{a+b}$, $v_0\,v_b$, $v_a\,v_{a+b}$. 
The outer vertices $u_0$, $u_a$, $u_b$, $u_{a+b}$ appear in $C$ in the sequence $u_0$, $u_a$, $u_{a+b}$, $u_b$, or  $u_0$, $u_a$,  $u_b$, $u_{a+b}$.
The edge $v_0\,v_b$ is placed in one of the subpaths of $C$ we obtain by removing the edges $u_0\,u_a$, $u_b\,u_{a+b}$; the same holds for the edge $v_a\,v_{a+b}$, and it may or may not belong to the same subpath as  $v_0\,v_b$. Together there are, up to symmetry, 48 different orderings of the vertices $u_0,u_a,u_b,u_{a+b},v_0,v_b ,v_a,v_{a+b}$ on $C$.

We show how the $4$-hooked construction works in the hypothesis that the vertices $u_0$, $u_a$, $u_b$, $u_{a+b}$ are ordered in $C$ in the sequence $u_0$, $u_a$, $u_{a+b}$, $u_b$, and that $v_0\,v_b$ belongs to the subpath $u_0\,P\,u_b$, whereas $v_a\,v_{a+b}$ is in $u_{a+b}\,P\,u_a$ ($u_0\,P\,u_b$,   $u_{a+b}\,P\,u_a$  are the subpaths of $C$ we obtain by removing the edges $u_0\,u_a$, $u_b\,u_{a+b}$). We also assume that $v_0$ precedes $v_b$ in the path $u_0\,P\,u_b$, and $v_{a+b}$ precedes $v_{a}$ in the path $u_{a+b}\,P\,u_{a}$.  Then,  by removing the edges  $v_0\,v_b$,
$v_a\,v_{a+b}$ in $C-\{u_0\,u_a, u_b\,u_{a+b}\}$, we obtain the following four subpaths: $v_0\,P\,u_0$, $v_b\,P\,u_b$, $v_{a+b}\,P\,u_{a+b}$, $v_{a}\,P\,u_{a}$.
We will also consider the subpaths $v_0\,P\,v_a$, $v_{b}\,P\,v_{a+b}$ we obtain from $C$ by removing the edges $v_0\,v_b$,  $v_a\,v_{a+b}$.

In $H_i$, with $1\leq i\leq\lambda-1$, we consider the subpaths $v^i_j\,P\,u^i_j$,  $j\in\{0, b, a, a+b\}$, which correspond to the above four subpaths of $C$.  In $H_{\lambda}$, we consider the subpaths $v^{\lambda}_0\,P\,v^{\lambda}_a$, $v^{\lambda}_{b}\,P\,v^{\lambda}_{a+b}$,
which corresponds to the subpaths $v_0\,P\,v_a$, $v_{b}\,P\,v_{a+b}$ of $C$.

For $1\leq i\leq\lambda-2$ and $j \in\{0, b, a, a+b\}$, we join the path $v^i_j\,P\,u^i_j$ in $H_i$ to the path $v^{i+1}_j\,P\,u^{i+1}_j$ in $H_{i+1}$ by the edge $u^i_j\,v^{i+1}_j$, 
and get a path $v^1_j\,P\,u^{\lambda-1}_j$ from $v^1_j$ in $H_1$ to $u^{\lambda-1}_j$ in $H_{\lambda-1}$. The union of the four paths is a disconnected graph covering all the vertices of $G$, with the exception for the vertices in $H_0\cup H_{\lambda}$.  We connect the four paths $v^1_j\,P\,u^{\lambda-1}_j$,  $j \in\{0, b, a, a+b\}$ to the paths $u^0_0\,P\,u^0_b$, $u^0_{a+b}\,P\,u^0_a$ in $H_0$ and to the paths $v^{\lambda}_0\,P\,v^{\lambda}_a$, $v^{\lambda}_{b}\,P\,v^{\lambda}_{a+b}$ in $H_{\lambda}$ by adding the spokes $u^0_j\,v^1_j$, $u^{\lambda-1}_j\,v^{\lambda}_j$, $j \in \{0, b, a, a+b\}$.
We thus obtain a Hamilton cycle in $G$, and the assertion follows. We summarize the construction in the diagram in Figure \ref{fig_4hooked_construction}. 

\begin{figure}
\begin{center}
\includegraphics[width=13cm]{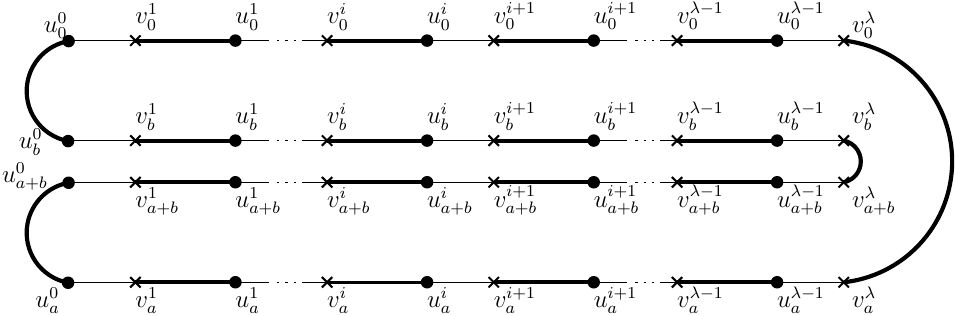}
\end{center}
\caption{The $4$-hooked construction for the generalized rose window graphs described in Proposition \ref{pro_4hooked}. The bold lines represent the paths $u^0_0\,P\,u^0_b$, $u^0_{a+b}\,P\,u^0_a$ in $H_0$,
$v^{\lambda}_0\,P\,v^{\lambda}_a$, $v^{\lambda}_{b}\,P\,v^{\lambda}_{a+b}$ in $H_{\lambda}$, and the paths 
$u^i_j\,P\,v^i_j$ in $H_i$, with $j \in\{0, b, a, a+b\}$ for $1\leq i\leq\lambda-1$; the paths are connected by the spokes $u^i_j\,v^{i+1}_j$, for $0\leq i\leq\lambda-1$.}
\label{fig_4hooked_construction} \end{figure}

We can repeat the same construction even if the four edges
$u_0\,u_a$, $u_b\,u_{a+b}$, $v_0\,v_b$, $v_a\,v_{a+b}$
are arranged on $C$ in a different way from that considered above when the outer and inner edges among these edges alternate on $C$, with the exceptions of the following  four orderings: 
$u_0, u_a, v_0, v_b, u_b, u_{a+b}, v_a, v_{a+b}$,
$u_0, u_a, v_b, v_0, u_b, u_{a+b}, v_{a+b}, v_{a}$,
$u_0, u_a, v_{a+b}, v_a, u_b, u_{a+b},$ $ v_b, v_{0}$,
$u_0, u_a, v_a, v_{a+b}, u_b, u_{a+b}, v_{0}, v_{b}$.
In such exceptions we find either a Hamilton path $v_0\,P\,v_a$ or a Hamilton path $u_0\,P\,u_b$ in $H_0$. For instance,  in the case of the sequence $u_0$, $u_a$,  $v_{0}$, $v_{b}$, $u_b$, $u_{a+b}$, $v_a$, $v_{a+b}$, we find a Hamilton path from $u_0$ to $u_b$ by removing the edges  $u_0\,u_a$, $u_b\,u_{a+b}$,  $v_a\,v_{a+b}$
and adding the edges $u_a\,v_a$, $u_{a+b}\,v_{a+b}$. The assertion then follows from Proposition \ref{pro_2hooked}.

When the outer and the inner edges from $\{u_0\,u_a$, $u_b\,u_{a+b}$, $v_0\,v_b$, $v_a\,v_{a+b}\}$ do not alternate on $C$, we can almost always find either a Hamilton path $v_0\,P\,v_a$ or a Hamilton path $u_0\,P\,u_b$ in $H_0$ (sometimes by relabeling the vertices appropriately) and then use Proposition \ref{pro_2hooked} to find a Hamilton cycle in the graph $G$.
\medskip

We cannot find such paths directly only if the cycle $C$ is elusive.
By Remark \ref{remark1_2} we may assume that the cycle $C$ is elusive of type 1.
Then by Lemma \ref{lem_4hookedmain} we can use the $2$-hooked construction or the $4$-hooked construction described above to obtain a Hamilton cycle in the graph $G$.
In the special case where $b\equiv -2a\pmod m$  or $a\equiv -2b\pmod m$, and the Hamilton cycle $C$ contains the subpaths $(v_b, v_0, u_0, u_a)$ 
and $(u_b, u_{a+b}, v_{a+b}, v_a, v_{a-b}, u_{a-b})$ 
occurring in this order in $C$, we may use the $2$-hooked construction as described in Remark \ref{rem_2hooked} and Example \ref{example}.  
\end{proof}

By combining Lemmas \ref{lem_classification_ham_cycle_Igraph}, 
 \ref{lemma_labeling} and Propositions \ref{pro_alternating_construction},  \ref{pro_2hooked},  \ref{pro_4hooked},  we can prove our main result that every  connected generalized rose window graph is hamiltonian.    
\medskip

\emph{Proof of Theorem \ref{th_rosewin_hamiltonian}.} 
Let $G=R(m;a,b,c)$ be a connected generalized rose window graph and let $H$ be the graph obtained from $G$ by removing the spokes of type $c$. Let $\lambda=\gcd(m, a, b)-1$.

First, we consider the case where $\lambda=0$, that is,  $H$ is a connected spanning subgraph of $G$. By Theorem \ref{th_cubicbicirculants}, we know that a connected $I$-graph is hamiltonian, with the exception of the generalized Petersen graphs $G(n, 2)$, with $n\equiv 5\pmod 6$. Therefore, if  $H$ is  not isomorphic to a graph $G(n, 2)$, then a Hamilton cycle of $H$ is also a Hamilton cycle of $G$. We find a Hamilton cycle in $G$ even if $H$ is a generalized Petersen graph $G(n, 2)$: Theorem \ref{thm:GPhamilton} assures the existence of a Hamilton path in $H$ connecting the vertices $u_0$ and $v_c$, which are adjacent in $G$ but not in $H$, since $H$ contains no spokes of type $c$; adding the spoke $u_0v_c$ yields a Hamilton cycle in $G$. Thus, the assertion follows if $\lambda=0$.
In the rest of the proof we consider $\lambda>0$.  

Let $H_0$ be the connected component of $H$ containing the vertex $u_0$.
If $H_0$ is not isomorphic to a generalized Petersen graph $G(n, 2)$, with $n\equiv 5\pmod 6$, then the assertion follows from 
Lemmas \ref{lem_classification_ham_cycle_Igraph},  
 \ref{lemma_labeling} and  Propositions \ref{pro_alternating_construction},  \ref{pro_2hooked},  \ref{pro_4hooked}.

Let us now consider the case where $H_0$ is the generalized Petersen graph $G(n, 2)$, with $n\equiv 5\pmod 6$.  Notice that $n=m/\gcd(m, a,  b)$ and that the indices of the vertices of $H_0$ are all multiples of $\gcd(m,a,b)=\lambda+1$. By Theorem \ref{thm:GPhamilton}, we can find a Hamilton path $v_0\,P\,u_x$  in $H_0$  connecting $v_0$ to $u_x$, and also a Hamilton path $v_x\,P\,u_0$, for every integer $x\in (\lambda+1)\mathbb Z_m$,  $x\not\equiv 0\pmod{m}$. 
For odd values of $\lambda$, we select $x\in (\lambda+1)\mathbb Z_m$ such that $x+(\lambda+1)c\not\equiv 0\pmod{m}$. For even values of $\lambda$, we do not add additional conditions on $x$, but we select another integer $y\in (\lambda+1) \mathbb Z_m$ such that $x\not\equiv y\pmod{m}$, $y+(\lambda+1)c\not\equiv 0\pmod{m}$, and consider the vertices $u^{\lambda}_y=u_{y+\lambda\,c}\in H_{\lambda}$, $v_{y+(\lambda+1)c}\in H_0$. Notice that the choice of $x$ and $y$ is always possible,  since $n\geq 5$.

We now construct a Hamilton cycle in $G$ as follows.  We take the path $v^i_{0}\,P\,u^i_{x}$ in $H_i$ with $i$ odd, $1\leq i\leq\lambda-1$,  and the path $v^i_{x}\,P\,u^i_{0}$ in each $H_i$ with $i$ even, $1\leq i\leq \lambda-1$. 
We join the paths by the spokes $u^i_{x}\,v^{i+1}_{x}$ for $1\leq i\leq \lambda-2$ with $i$ odd, and $u^i_{0}\,v^{i+1}_{0}$ for $1\leq i\leq \lambda-2$ with $i$ even. 
For odd values of $\lambda$, we obtain a path $v^1_{0}\,P\,u^{\lambda-1}_{0}$ connecting the vertices $v^1_0$ and  $u^{\lambda-1}_0$; for even values of $\lambda$, we have a path
$v^1_{0}\,P\,u^{\lambda-1}_{x}$ connecting the vertices $v^1_0$ and  $u^{\lambda-1}_{x}$; 
Both paths  $v^1_{0}\,P\,u^{\lambda-1}_{0}$ and $v^1_{0}\,P\,u^{\lambda-1}_{x}$ cover all the vertices in $G-(H_0\cup H_{\lambda})$.

For odd values of $\lambda$, we take the Hamilton path $v^{\lambda}_0\,P\,u^{\lambda}_x$  in $H_{\lambda}$, 
and the Hamilton path $v^0_{x+(\lambda+1)c}\,P\,u^0_0$ in $H_0$  (whose existence follows from Theorem \ref{thm:GPhamilton} by the assumptions on $x$).  We join the paths $v^0_{x+(\lambda+1)c}\,P\,u^0_0$, 
$v^1_{0}\,P\,u^{\lambda-1}_0$, $v^{\lambda}_0\,P\,u^{\lambda}_x$ by the spokes $u^0_0\,v^1_0$, 
$u^{\lambda-1}_0\,v^{\lambda}_0$, $u^{\lambda}_x\,v^0_{x+(\lambda+1)c}$, and obtain
a Hamilton cycle in $G$.

For even values of $\lambda$, we take the Hamilton path $v^{\lambda}_x\,P\,u^{\lambda}_y$  in $H_{\lambda}$
and the Hamilton path $v^0_{y+(\lambda+1)c}\,P\,u^0_0$ in $H_0$  (whose existence follows from Theorem \ref{thm:GPhamilton} by the assumptions on $y$).
We join the paths $v^0_{y+(\lambda+1)c}\,P\,u^0_0$,  $v^1_0\,P\,u^{\lambda-1}_x$,  $v^{\lambda}_x\,P\,u^{\lambda}_y$
by the spokes $u_0\,v^1_0$, $u^{\lambda-1}_x\,v^{\lambda}_x$, $u^{\lambda}_y\,v^0_{y+(\lambda+1)c}$, and obtain
a Hamilton cycle in $G$, which completes the proof.  We summarize the construction in the diagram in Figure \ref{fig_generalized_petersen}.
\qed

\section{Concluding remarks}

Proving that all generalized rose window graphs are hamiltonian could be the first step to proving Conjecture \ref{conjectB}.
The next step would be to consider the pentavalent generalized  Taba\v cjn graphs, which are obtained from the generalized rose window graphs by adding an additional set of spokes, similarly as the generalized  rose window graphs are obtained from the $I$-graphs, see \cite{Taba1,Taba2}. Given $m \ge 3$ and $a,b,c,d \in \ZZ_m \setminus\{0\}$ with $a,b \ne m/2$,
the \emph{generalized Taba\v cjn graph} $T(m;a,b,c,d)$ is defined to be the bicirculant graph $B(m;\{ a,-a\},\{0,c,d\},\{b,-b\})$. 
Every connected generalized Taba\v cjn graph $T(m;a,b,c,d)$ contains three generalized rose window graphs as subgraphs, namely $R(m;a,b,c)$, $R(m;a,b,d)$ and $B(m;R,S \setminus \{0\},T)$, which is isomorphic to $R(m;a,b, d-c)$ by Proposition \ref{pro:subtract_isomorphism}. 
It may happen that at least one of these is connected. In this case, also the graph $T(m;a,b,c,d)$ is hamiltonian by Theorem \ref{th_rosewin_hamiltonian}. 
Moreover, the graph $T(m;a,b,c,d)$ contains the cubic Haar graph  $B(m;\emptyset,\{0,c,d\} ,\emptyset)$ as a subgraph. If that graph is connected, it is hamiltonian by \cite[Theorem 3.1]{AlZh1989}; this happens when $\gcd(m,c,d )=1$. 

We can apply the same reasoning to more general bicirculant graphs: if a bicirculant contains a connected generalized rose window graph as a subgraph, then it is hamiltonian by Theorem \ref{th_rosewin_hamiltonian}; if it contains a connected cubic Haar graph as a subgraph, then it is hamiltonian by \cite[Theorem 3.1]{AlZh1989}.
 
\begin{proposition} \label{prop:primesHaar}
Let $G=H(m;S)$ be a connected cyclic Haar graph with $|S| \ge 4$. If  $m$ is a product of at most three prime powers, then $G$ is hamiltonian. 
\end{proposition}
\begin{proof}
Let $S=\{0,c_1,\dots, c_{s-1}\}$, where $s=|S| \ge 4$. Since the graph $G$ is connected, we have $\gcd(m,S)=1$. 
If   $G$  contains a connected  cubic Haar graph as a subgraph, then it is hamiltonian. 
Therefore we assume that $G$ does not contain  a connected cubic Haar graph as a subgraph and we will show that in this case $m$ needs to be a product of at least four prime powers. 

Since $\gcd(m,c_1,c_2) > 1$, it is divisible by some prime, say $p$. Since $\gcd(m,S)=1$, there exists an element of $S$, say $c_i$, that is not divisible by $p$. Therefore there exists another prime, say $q$, such that $\gcd(m,c_1,c_i)$ is divisible by $q$. Now there exists some element of $S$ that is not divisible by $q$, say $c_j$ (it may happen that $c_j=c_2$). Therefore there exists another prime, say $r$, such that $\gcd(m,c_i,c_j)$ is divisible by $r$. Thus $m$ is a product of at least three prime powers. 

Suppose that $m$ is a product of exactly three prime powers, namely, the powers of $p$, $q$ and $r$. Since $\gcd(m,S)=1$, again there exists an element of $S$, say $c_k$, that is not divisible by $r$ (it may happen that $c_j=c_1$ or $c_j=c_2$). Now we have elements $c_i,c_j,c_k$ from $S$ such that $c_i$ is not divisible by $p$, $c_j$ is not divisible by $q$ and $c_k$ is not divisible by $r$. 
On the other hand all of $\gcd(m,c_i,c_j)$, $\gcd(m,c_i,c_k)$, $\gcd(m,c_j,c_k)$ are greater than one. That means that $c_i,c_j$ are both divisible by $r$, $c_i,c_k$ are both divisible by $q$ and $c_j,c_k$ are both divisible by $p$.  But then $c_i-c_k$ is  not divisible by any of $p,r$ and $c_j-c_k$ is not divisible by $q$. It follows that $\gcd(m, c_i-c_k,c_j-c_k)=1$ and $G$ contains a connected cubic Haar graph $H(m;\{0,c_i-c_k,c_j-c_k\})$ as a subgraph, a contradiction. Therefore $m$ is a product of at least four prime powers.
\end{proof}

\begin{proposition} \label{prop:primes}
Let $G=B(m;a,S,b)$ be a connected bicirculant with $|S| \ge 3$ and $a,b \ne m/2$. If  $m$ is a product of at most three prime powers, then $G$ is hamiltonian. 
\end{proposition}
\begin{proof}
Let $S=\{0,c_1,\dots, c_{s-1}\}$, where $s=|S| \ge 3$. Since the graph $G$ is connected, we have $\gcd(m,a,S,b)=1$. If the graph  $G$ contains a connected generalized rose window graph or a connected cubic Haar graph as a subgraph, then it is hamiltonian. Therefore we assume that this is not the case and we will show that then $m$ needs to be a product of at least four prime powers. 

We may assume that $\gcd(m,S)= d >1$, otherwise already $B(m; \emptyset, S, \emptyset)$ is connected and the claim follows from Proposition \ref{prop:primesHaar}.
Therefore there exists a prime $r$ that divides $d$. Since the graph $G$ is connected and it does not contain a connected generalized rose window graph as a subgraph, at least one of $a,b$, say $a$, must be coprime to $r$; therefore there exists  a prime  $p$ that is coprime to $d$ such that $p$ is coprime to $r$ and $p$ divides $\gcd(m,a,b,c_1)$. In particular $p \ne r$. Since $p$ is coprime to $d$, there exists $c_i \in S \smallsetminus \{c_1\}$ that is not divisible by $p$. Therefore there exists a third prime, say $q$, that divides $\gcd(m,a,b,c_i)$. Thus $m$ is a product of at least three prime powers. 

Suppose that $m$ is a product of exactly three prime powers, namely, the powers of $p$, $q$ and $r$. Since the graph $G$ is connected, there exists $c_j \in S \smallsetminus \{c_i\}$ that is not divisible by $q$ (it may happen that $c_j=c_1$). Since $\gcd(m;a,b,c_j)>1$, it follows that $c_j$ must be divisible by $p$. Now we have elements $c_i,c_j$ from $S$ such that $c_i$ is divisible by $q$ and is coprime to $p$, $c_j$ is divisible by $p$ and is coprime to $q$. But then $c_i-c_j$ is  not divisible by any of $p,q$ and $a$ is not divisible by $r$. It follows that $\gcd(m, a,c_i-c_j,b)=1$ and $G$ contains a connected rose window graph $S(m;a,\{0,c_i-c_j\},b)$ as a subgraph, a contradiction. Therefore $m$ is a product of at least four prime powers.
\end{proof}

\section*{Acknowledgements}

The authors would like to thank Brian Alspach for a careful reading of the manuscript and for suggesting a number of grammatical corrections and improvements in the presentation of the results.
 
Simona Bonvicini is a member of GNSAGA of Istituto Nazionale di Alta Matematica (INdAM).
Toma\v{z} Pisanski is supported in part by the Slovenian Research Agency (research program P1-0294 and research projects J1-4351, J5-4596 and BI-HR/23-24-012).
Arjana \v Zitnik is supported in part by the Slovenian Research Agency (research program P1-0294 and research projects J1-3002 and J1-4351).

\section*{Conflict of interest}
The authors declare that they have no conflict of interest.



\appendix

\section{Proof of Lemma \ref{lem_4hookedmain}}
\label{appendix}

In this Section we prove Lemma \ref{lem_4hookedmain}. 
We have the following assumptions.
\begin{assumption} \label{assumption1}
Given an $I$-graph $G=I(m,a,b)$,  assume that  
\begin{itemize}
\item $a \ne\pm b $, 
\item the graph  $G$  contains a $4$-hooked Hamilton cycle $C$ that is elusive and the vertices $u_0$, $u_a$, $u_b$, $u_{a+b}$, $v_{a+b}$, $v_a$, $v_b$, $v_0$ are ordered on $C$   in this way. 
\end{itemize}
\end{assumption}
Note that because of the ordering of the vertices in the elusive Hamilton cycle,  the vertices $u_b$ and $v_b$ are not adjacent in $C$ and also $u_a$ and $v_a$ are not adjacent in $C$, therefore the vertex $u_{-a+b}$ is adjacent to the vertex $u_b$ in $C$, and $v_{a-b}$ is adjacent to $v_a$ in $C$.
\medskip

By Remark \ref{remark_lemma8}, it is enough to consider three cases: 
\begin{enumerate}[I.]
\item   the path  $(u_a,u_{0}, u_{-a})$ is contained in $C$ and $u_{-a}$ is not adjacent to $v_{-a}$ on $C$, i.e.
the path $(v_{-a+b}, v_{-a}, v_{-a-b})$ is contained in $C$,  
\item the path $(u_a,u_{0}, u_{-a}, v_{-a})$ is contained in $C$,
\item   the cycle $C$ contains the subpaths 
$(v_b, v_0, u_0, u_a,u_{2a},v_{2a})$ 
and $(u_b, u_{a+b},$ $ v_{a+b},$ $v_a,$ $v_{a-b}, u_{a-b})$
occurring in this order in $C$.
\end{enumerate}

We deal with case I in Lemma \ref{lemma_special1}, with case II in Lemma \ref{lemma_special2} and with case III  in Lemma \ref{lemma_special3b}.

In order to find the appropriate paths or cycles, we will sometimes need to consider additional vertices on the cycle $C$,  such as $u_{\pm 2a}$, $v_{\pm 2a}$, $u_{-a+b}$, $v_{-a+b}$ and others. 
In some cases the  arrangement of the vertices $u_0$, $u_a$, $u_b$, $u_{a+b}$, $v_{a+b}$, $v_a$, $v_b$, $v_0$ on $C$ forces these vertices to be different from others. 
It is possible that for some $m$ some of the vertices we will consider later are not necessarily distinct from the existing ones --  for example  $v_a=v_{-a+b}$ when $b=2a$ -- however our construction is still valid.
We denote with $u_a\,P\,u_b$, $u_{a+b}\,P\,v_{a+b}$, $v_a\,P\,v_b$, $v_0\,P\,u_0$ the paths we obtain from $C$ by removing the edges $u_0u_a$, $u_bu_{a+b}$, $v_0v_b$ and $v_av_{a+b}$. 

\smallskip

\begin{lemma}  \label{lemma_special1}
Given an $I$-graph $G=I(m,a,b)$, in addition to Assumption \ref{assumption1}, assume that the path $(v_{-a+b}, v_{-a}, v_{-a-b})$ is contained in $C$. 
Then the graph $G$ contains a standard $4$-hooked Hamilton cycle or a $2$-hooked Hamilton cycle.
\end{lemma}
\begin{proof}
When the path 
$(v_{-a+b}, v_{-a}, v_{-a-b})$ is in $C$, the  edge 
$u_{-a}v_{-a}$ is not in $C$, since the vertex $v_{-a}$ already has two neighbours in $C$. Therefore the edge $u_0u_{-a}$ is in $C$, which in turn means that the edge $u_0v_0$ is not in $C$. Thus, also the edge $v_0v_{-b}$ is in $C$. 
The path $(v_{-a+b}, v_{-a}, v_{-a-b})$ is contained in one of the four paths $u_a\,P\,u_b$, $u_{a+b}\,P\,v_{a+b}$, $v_a\,P\,v_b$, $v_0\,P\,u_0$ in two different orientations. 
We deal with each of the four situations separately. \\[-3mm]

\textbf{Case (a).} If $(v_{-a+b}, v_{-a}, v_{-a-b})$ is in $u_a\,P\,u_b$, then the vertices occur in the order  
$u_a$, $v_{-a+b}$, $v_{-a}$, $v_{-a-b}$, $u_b$, or $u_a$, $v_{-a-b}$, $v_{-a}$, $v_{-a+b}$, $u_b$. In the former case we remove the edges $v_{-a}v_{-a-b}$,
$u_0u_{-a}$, $v_0v_{-b}$, add the chords $u_{-a}v_{-a}$, $u_0v_0$, and obtain a Hamilton path from $v_{-b}$ to $v_{-a-b}$ yielding a Hamilton path from $v_0$ to $v_a$. In the latter case, we obtain a new labeling of $C$ by adding $a$ modulo $m$ to the subscripts of the vertices of $C$; the edges $u_0u_a$, $v_0v_b$, $u_bu_{a+b}$, $v_av_{a+b}$ are still in $C$, and the vertices occur in the order $u_0$, $u_a$, $v_0$, $v_b$, $u_b$, $u_{a+b}$, $v_{a+b}$, $v_a$, so the cycle $C$ is standard.

\textbf{Case (b).} If $(v_{-a+b}, v_{-a}, v_{-a-b})$ is in $u_{a+b}\,P\,v_{a+b}$, then the vertices occur in the order  
$u_{a+b}$, $v_{-a+b}$, $v_{-a}$, $v_{-a-b}$, $v_{a+b}$, or $u_{a+b}$, $v_{-a-b}$, $v_{-a}$, $v_{-a+b}$, $v_{a+b}$. In the former case we remove the edges 
 $u_bu_{-a+b}$, $v_0v_{b}$, $v_{-a}v_{-a+b}$, add the chords $u_bv_b$, $u_{-a+b}v_{-a+b}$, and obtain a Hamilton path from $v_0$ to $v_{-a}$ providing a Hamilton path from $v_0$ to $v_a$. In the latter case we remove the edges $u_bu_{-a+b}$, $u_0u_{-a}$,
$v_{-a}v_{-a+b}$, add the edges $u_{-a}v_{-a}$, $u_{-a+b}v_{-a+b}$, and find a Hamilton path from $u_0$ to $u_{b}$. \\[-3mm]

\textbf{Case (c).} If $(v_{-a+b}, v_{-a}, v_{-a-b})$ is in $v_{a}\,P\,v_{b}$, then the vertices occur in the order  
$v_{a}$, $v_{-a+b}$, $v_{-a}$, $v_{-a-b}$, $v_{b}$, or $v_{a}$, $v_{-a-b}$, $v_{-a}$, $v_{-a+b}$, $v_{b}$. In the former case we remove the edges $u_bu_{-a+b}$, $v_0v_{b}$, $v_{-a}v_{-a+b}$, add the edges $u_bv_b$, $u_{-a+b}v_{-a+b}$, and obtain a Hamilton path from
$v_0$ to $v_{-a}$ yielding a Hamilton path from $v_0$ to $v_a$. In the latter case we remove the edges $u_0u_{-a}$, $v_0v_{-b}$, $v_{-a}v_{-a-b}$, add the chords $u_{-a}v_{-a}$, $u_0v_0$, and find a Hamilton path from $v_{-b}$ to $v_{-a-b}$ that provides a Hamilton path from $v_0$ to $v_a$. 

\textbf{Case (d).}
If $(v_{-a+b}, v_{-a}, v_{-a-b})$ is in $v_{0}\,P\,u_{0}$ and the vertices occur in the sequence 
$v_0$, $v_{-a-b}$, $v_{-a}$, $v_{-a+b}$, $u_{0}$, then we remove the edges $u_0u_{-a}$, $u_{b}u_{-a+b}$, $v_{-a}v_{-a+b}$, add the chords $u_{-a}v_{-a}$, $u_{-a+b}v_{-a+b}$, and find a Hamilton path from $u_0$ to $u_b$. 

The case where $(v_{-a+b}, v_{-a}, v_{-a-b})$ is in $v_{0}\,P\,u_{0}$ and the vertices appear in the sequence $v_0$, $v_{-a+b}$, $v_{-a}$, $v_{-a-b}$, $u_{0}$ deserves special attention; it depends on the presence of the edge $u_{-b}v_{-b}$ in $C$. 
More specifically, if the edge $u_{-b}v_{-b}$ is in $C$, then we find one of the two paths $(u_{a-b}, u_{-b}, v_{-b})$ or $(u_{-a-b}, u_{-b}, v_{-b})$ in $C$. 

If $C$ contains the path $(u_{a-b}, u_{-b}, v_{-b})$, then $u_{-a-b}$ is adjacent to $v_{-a-b}$ in $C$. We can thus remove the edges $u_{-a-b}v_{-a-b}$, $u_{-b}v_{-b}$, add the chord $u_{-b}u_{-a-b}$, and find a Hamilton path from $v_{-b}$ to $v_{-a-b}$ that yields a Hamilton path from $v_0$ to $v_a$.

If $(u_{-a-b}, u_{-b}, v_{-b})$ is in $C$, then $u_{a-b}$ is adjacent to $v_{a-b}$ in $C$ (which is adjacent to $v_a$). In this setting we remove the edges $u_{-b}v_{-b}$, $u_{a-b}v_{a-b}$, add the edge $u_{-b}u_{a-b}$, and obtain a Hamilton path from $v_{-b}$ to $v_{a-b}$ providing a Hamilton path from $v_0$ to $v_a$.

If the edge $u_{-b}v_{-b}$ is not in $C$, then the path $(u_{-a-b}, u_{-b}, u_{a-b})$ is in $C$, and can be contained in one of the four paths 
$u_a\,P\,u_b$, $u_{a+b}\,P\,v_{a+b}$, $v_a\,P\,v_b$, $v_0\,P\,u_0$. 
In each case we can find a Hamilton path from $u_0$ to $u_b$, or from $v_0$ to $v_a$, or define a new labeling of the vertices of $C$ by adding a suitable integer modulo $m$ to the subscripts of the vertices, which 
makes the cycle $C$ standard.
This is the case where $(u_{-a-b}, u_{-b}, u_{a-b})$ is in $v_a\,P\,v_b$ and the vertices appear in the sequence $v_a$, $u_{a-b}$, $u_{-b}$, $u_{-a-b}$, $v_b$: we add $b$ modulo $m$ to the subscripts of the vertices in $C$ and obtain a new labeling of $C$ which still gives the edges $u_0u_a$, $v_0v_b$, $u_bu_{a+b}$, $v_av_{a+b}$ in $C$, and the vertices occur in the order $u_0$, $u_a$, $v_a$, $v_{a+b}$, $u_{a+b}$, $u_b$, $v_0$, $v_b$. For the other cases we summarize the construction in Figure \ref{fig_4hooked_position_u-b}.  
We could also notice that by adding   $b$ modulo $m$ to the subscripts of the vertices in Figures \ref{fig_4hooked_position_u-b}  (a)-(f), the new labeling makes the cycle $C$  standard. However, in Figures \ref{fig_4hooked_position_u-b} (g)-(h) the new labeling keeps $C$ elusive. 
\begin{figure}[htb]
\begin{center}
\includegraphics[width=16cm]{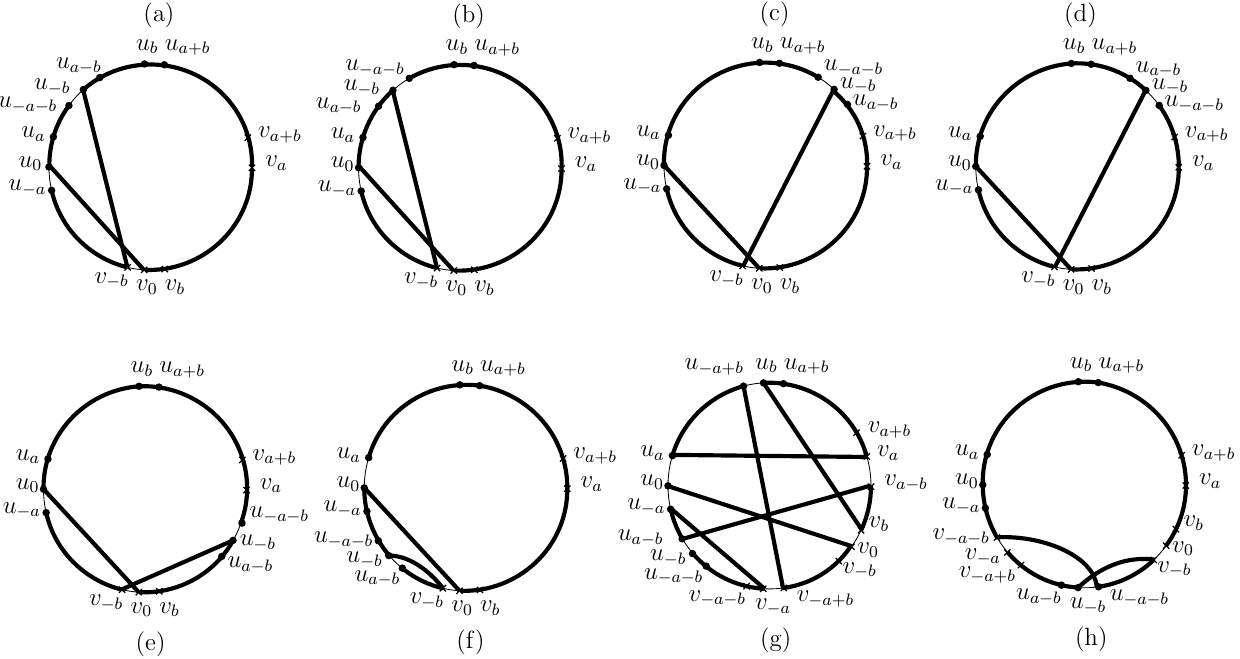}
\end{center}
\caption{The bold lines denote a Hamilton path, connecting two outer or two inner vertices, that provides a Hamilton path from $u_0$ to $u_b$ or $v_0$ to $v_a$ in the proof of Lemma  \ref{lemma_special1}
when the path $(u_{-a-b}, u_{-b}, u_{a-b})$ is in $C$. The Hamilton paths from $u_0$ to $u_b$ or from $v_0$ to $v_a$ can be obtained  by adding a suitable integer modulo $m$ to the subscripts of the vertices; the path $(u_{-a-b}, u_{-b}, u_{a-b})$ can be contained in one of the paths $u_a\,P\,u_b$, $u_{a+b}\,P\,v_{a+b}$, $v_a\,P\,v_b$, $v_0\,P\,u_0$, and its vertices can occur in two different orders. The figure does not show the case where $(u_{-a-b}, u_{-b}, u_{a-b})$ is contained in the path $v_a\,P\,v_b$ and the vertices occur in the order $v_a$, $u_{a-b}$, $u_{-b}$, $u_{-a-b}$, $v_{b}$ because from this case we can obtain a new labeling of the vertices of $C$ that makes the cycle $C$ standard.
}\label{fig_4hooked_position_u-b}
\end{figure}
\end{proof}

\begin{lemma}  \label{lemma_special2}
Given an $I$-graph $G=I(m,a,b)$, in addition to Assumption \ref{assumption1}, assume that the path $(u_{0}, u_{-a}, v_{-a})$ is contained in $C$.
Then the graph $G$ contains a standard $4$-hooked Hamilton cycle or a $2$-hooked Hamilton cycle.
\end{lemma}
\begin{proof}
When the path $(u_0, u_{-a}, v_{-a})$ is in the cycle $C$, the edge $u_0v_0$ is not in  $C$, so the edge $v_0v_{-b}$ is in $C$. The analysis of this case again depends on the presence of the edge $u_{-b}v_{-b}$ in $C$.

\textbf{Case (a).} Assume that the edge $u_{-b}v_{-b}$ is in $C$. If the edge $u_{-b}u_{a-b}$ is not in $C$, then $C$ contains the edge $u_{a-b}v_{a-b}$, and we can find a Hamilton path from $v_{-b}$ to $v_{a-b}$ by deleting the edges $u_xv_x$ with $x\in\{a-b, -b\}$, and by adding the edge  $u_{-b}u_{a-b}$; the path yields a Hamilton path from $v_0$ to $v_a$. We consider the case where $u_{-b}u_{a-b}$ is in $C$. Then $C$ also contains the edge $u_{-a-b}v_{-a-b}$, and might contain $v_{-a}v_{-a+b}$ or not. 

If $v_{-a}v_{-a+b}$ is not in $C$, then the edge $u_{-a+b}v_{-a+b}$ is in $C$, 
and we can find a Hamilton path from $u_{-a}$ to $u_{-a+b}$ by removing the edges $u_xv_x$ with $x\in\{-a+b, -a\}$ from $C$, and by adding the edge $v_{-a}v_{-a+b}$; the path yields a Hamilton path from $u_0$ to $u_b$. Consider the case where $v_{-a}v_{-a+b}$ is in $C$ (so $v_{-a-b}u_{-a-b}$ is in $C$). We also assume that the vertices $v_{-a}$,
$u_{-a}$, $u_{-a-b}$, $v_{-a-b}$ occur in $C$ in this order, otherwise we find a Hamilton path from $u_{-a}$ to $u_{-a-b}$ providing a Hamilton path from $u_0$ to $u_b$ (we remove the edges $u_xv_x$ with $x\in\{-a-b, -a\}$, and add the edge $v_{-a}v_{-a-b}$). 
Whatever the position of edge $u_{-a-b}v_{-a-b}$ in $C$, we remove the following edges from $C$: $u_bu_{-a+b}$,
$v_av_{a-b}$, $u_{-a-b}v_{-a-b}$, $v_0v_b$, $u_{-b}u_{a-b}$, $u_{-a}v_{-a}$, $u_0u_a$; we add the edges:
$u_0v_0$, $u_av_a$, $u_bv_b$, $u_{a-b}v_{a-b}$, $u_{-b}u_{-a-b}$, $v_{-a}v_{-a-b}$. We find a Hamilton path from $u_{-a}$ to $u_{-a+b}$,
which yields a Hamilton path from $u_0$ to $u_b$, if the edge $u_{-a-b}v_{-a-b}$ is not in $v_0\,P\,u_0$. 

If the edge $u_{-a-b}v_{-a-b}$ is in $v_0\,P\,u_0$, then we add the edge $u_{-a}v_{-a}$ back, add also the edges $v_{-a+b}u_{-a+b}$,  $v_{-a}v_{-a-b}$ ,
remove the edge $v_{-a}v_{-a+b}$ and find a new Hamilton cycle, say $C'$, containing the subpath $(u_{-a-b}, u_{-b}, v_{-b}, v_0, u_0, u_{-a}, v_{-a}, v_{-a-b})$. We relabel the vertices of $C'$ by adding $(a+b)$ modulo $m$ to the subscripts of the vertices; 
we obtain a new labeling of $C'$ having the edges $u_0u_a$, $v_0v_b$, $u_bu_{a+b}$, $v_av_{a+b}$, and the vertices occurring in the order 
$u_0$, $u_a$, $v_a$, $v_{a+b}$, $u_{a+b}$, $u_{b}$, $v_{b}$, $v_0$,
so the cycle $C'$ is standard.

\textbf{Case (b).} Assume that the edge $u_{-b}v_{-b}$ is not in $C$.
We note that the constructions in Figure \ref{fig_4hooked_position_u-b} (a)--(f) are independent from the position of the vertex $v_{-a}$ in $C$ ($v_{-a}$ does not appear in the figure). For this reason we can use the same constructions even in the case where $(u_0, u_{-a}, v_{-a})$ is in $C$. In order to conclude the proof of the case in which $(u_0, u_{-a}, v_{-a})$ is in $C$, it remains to consider the case where $(u_{-a-b}, u_{-b}, u_{a-b})$ is in $v_0\,P\,u_0$, and the vertices occur in the order $v_0$, $u_{-a-b}$, $u_{-b}$, $u_{a-b}$, $u_0$ (such a case replaces the constructions in (g)--(h) of Figure \ref{fig_4hooked_position_u-b}). 

If $(u_{-a}, v_{-a}, v_{-a+b})$ is in $v_0\,P\,u_0$, then $v_{-a-b}$ is adjacent to $u_{-a-b}$. In this setting we remove the edges $u_{-a}v_{-a}$, $u_{-a-b}v_{-a-b}$, add the chord $v_{-a}v_{-a-b}$, and find a Hamilton path from $u_{-a}$ to $u_{-a-b}$ that yields a Hamilton path from $u_0$ to $u_b$. 

If $(u_{-a}, v_{-a}, v_{-a-b})$ is in $v_0\,P\,u_0$, then we add $(a+b)$ modulo $m$ to the subscripts of the vertices and find a new labeling of $C$, which still gives the edges $u_0u_a$,
$v_0v_b$, $u_bu_{a+b}$, $v_av_{a+b}$ in $C$, and the vertices appear in the order $u_0$, $u_a$, $v_0$, $v_b$, $u_b$, $u_{a+b}$, $v_{a+b}$, $v_a$, so the cycle $C$ is standard.
This completes the proof.
\end{proof}

\begin{lemma}  \label{lemma_special3b}
Given an $I$-graph $G=I(m,a,b)$, in addition to Assumption \ref{assumption1}, assume that the paths 
$(v_b, v_0, u_0, u_a, u_{2a},v_{2a})$ and
$(u_b, u_{a+b},$ $ v_{a+b},$ $v_a,$ $v_{a-b}, u_{a-b})$
are contained in $C$. 
Then the graph $G$ contains a standard $4$-hooked Hamilton cycle or a $2$-hooked Hamilton cycle or $(b \equiv -2a\pmod m$ or $a \equiv -2b\pmod m)$  and  the cycle $C$  contains the subpaths $(v_b, v_0, u_0, u_a)$
and $(u_b, u_{a+b},$ $ v_{a+b}, v_a, v_{a-b}, u_{a-b}) $
occurring in this order in $C$.
\end{lemma}
\begin{proof}
If $(b \equiv -2a\pmod m$ or $a \equiv -2b\pmod m)$ then the assertion follows. Therefore we may assume that this is not the case.
This also  means that the vertices $u_{2a+b}$ and $u_0$ are distinct and the vertices  $v_{2a+b}$ and $v_0$ are distinct.

Since the vertices  $u_{a+b}$, $v_{a+b}$ are adjacent in $C$, also the vertices  $u_{2a+b},v_{2a+b}$ are adjacent in $C$. If the vertices occur in $C$ in order $u_{a+b}$, $v_{a+b}$, $u_{2a+b}$, $v_{2a+b}$, then we have a Hamilton path from $v_{a+b}$ to $v_{2a+b}$ yielding a Hamilton path from $v_0$ to $v_a$ by removing the edges $u_xv_x$ with $x\in\{a+b, 2a+b\}$ and adding the edge $u_{a+b}u_{2a+b}$.
In the following, we assume that the vertices occur in order $u_{a+b}$, $v_{a+b}$, $v_{2a+b}$, $u_{2a+b}$ and we distinguish two cases: 
(a) the subpath $(u_{2a}, v_{2a},v_{2a+b}, u_{2a+b})$ is in $C$ and
(b) the subpath $(u_{2a}, v_{2a},v_{2a-b})$ is in $C$.

\textbf{Case (a)} We assume that the subpath $(u_{2a}, v_{2a},v_{2a+b}, u_{2a+b})$ is in $C$. Then the edge $u_{2a-b}v_{2a-b}$ is in $C$.
In the case where the vertices $u_{2a}$, $v_{2a}$, $u_{2a-b}$, $v_{2a-b}$ occur in $C$ in this order, 
we find a Hamilton path from $u_{2a}$ to $u_{2a-b}$ yielding a Hamilton path from
$u_0$ to $u_b$ by removing the edges $u_xv_x$ with $x\in\{2a, 2a-b\}$  and adding the edge $v_{2a}v_{2a-b}$. 
In the case where the vertices occur in the order $u_{2a}$, $v_{2a}$, $v_{2a-b}$, $u_{2a-b}$,   we remove the edges
$u_xv_x$ with $x\in\{a-b, 2a-b\}$, add the edge $u_{a-b}u_{2a-b}$, and find a Hamilton path from $v_{a-b}$ to $v_{2a-b}$ 
providing a Hamilton path from $v_0$ to $v_a$.

\textbf{Case (b)} We assume that the subpath $(u_{2a}, v_{2a},v_{2a-b})$ is in $C$. If $(v_{a-b}, u_{a-b}, u_{-b}, v_{-b})$ is also in $C$, then $u_{2a-b}v_{2a-b}$ is in $C$. We remove the edges $u_xv_x$ with $x\in\{a-b, 2a-b\}$ from $C$, add the edge $u_{a-b}u_{2a-b}$,
and find a Hamilton path from $v_{a-b}$ to $v_{2a-b}$ that provides a Hamilton path from $v_0$ to $v_a$.
Finally, we consider the case where $(v_{a-b}, u_{a-b}, u_{2a-b})$ is in $C$.
If the vertices $v_0$, $u_0$, $v_{-b}$, $u_{-b}$ occur in this order in $C$, we  remove the edges $u_xv_x$ with $x\in\{0, -b\}$ from $C$, and add the edge $v_0v_{-b}$. The resulting Hamilton path from $u_0$ to $u_{-b}$ yields a Hamilton path from $u_0$ to $u_b$. 
If  the vertices $v_0$, $u_0$, $v_{-b}$, $u_{-b}$ occur in order  $v_0$, $u_0$, $u_{-b}$, $v_{-b}$ on $C$, 
the assertion follows from Lemma \ref{lemma_special4}. 
\end{proof}

For the remainder of this section we   now  have the following  assumptions.

\begin{assumption} \label{assumption2}
Given an $I$-graph $G=I(m,a,b)$, in addition to Assumption \ref{assumption1},  assume 
\begin{itemize}
\item $b \not \equiv -2a\pmod m$ and $a \not \equiv -2b\pmod m$,
\item the paths $(v_b, v_0, u_0, u_a, u_{2a},v_{2a})$
and $(u_b, u_{a+b},$ $ v_{a+b},$ $v_a,$ $v_{a-b}, u_{a-b})$ are   contained in $C$ in this order (so the edge $v_0 v_{-b}$ is not contained in $C$),
\item the edges $u_{a+b}v_{a+b}$ and $u_{2a+b}v_{2a+b}$ are contained in $C$,
\item the vertices $u_{a+b},v_{a+b},v_{2a+b}u_{2a+b}$ appear on $C$ in this order,
\item the paths 
$(u_{2a},v_{2a},v_{2a-b})$ and $(v_{a-b},u_{a-b},u_{2a-b})$ are contained in $C$, and
\item the vertices $v_0,u_0,u_{-b},v_{-b}$ appear on $C$ in this order.
\end{itemize}
\end{assumption}

\begin{lemma}  \label{lemma_special4}
Given an $I$-graph $G=I(m,a,b)$ containing a Hamilton cycle $C$ satisfying Assumption \ref{assumption2}, the graph $G$ contains a standard $4$-hooked Hamilton cycle or a $2$-hooked Hamilton cycle.
\end{lemma}

\begin{proof}
Note that in this setting, the edge $u_{-b}v_{-b}$ belongs to the cycle $C$. We distinguish three cases: (a) the edge $u_{-b}v_{-b}$ belongs to the subpath $u_a\,P\,u_b$ of the cycle $C$, (b) both edges $u_{-b}v_{-b}$ and $u_{2a+b}v_{2a+b}$ belong to the subpath $v_a\,P\,v_b$ and (c) the edge $u_{-b}v_{-b}$ belongs to the subpath $v_a\,P\,v_b$  and the edge $u_{2a+b}v_{2a+b}$ belongs to the subpath $u_a\,P\,u_b$.
\medskip

\textbf{Case (a).}  Assume that the edge $u_{-b}v_{-b}$ belongs to the subpath $u_a\,P\,u_b$ of the cycle $C$.
If  the vertices $v_{2a+b}$, $u_{2a+b}$ are also in $u_a\,P\,u_b$ and precede $u_{-b}$, $v_{-b}$, then we find a Hamilton path from $v_{a-b}$ to $v_{2a-b}$ providing a Hamilton path from $v_0$ to $v_a$ by
removing the edges $v_0v_b$, $v_{2a}v_{2a-b}$, $u_{2a+b}v_{2a+b}$, $u_{-b}v_{-b}$, $u_bu_{a+b}$, $v_{a-b}u_{a-b}$, and adding
the edges $v_{2a}v_{2a+b}$, $u_{a+b}u_{2a+b}$, $u_{-b}u_{a-b}$, $v_0v_{-b}$, $u_bv_b$. 
In the case where $v_{2a+b}$, $u_{2a+b}$ come after $u_{-b}$, $v_{-b}$ in $u_a\,P\,u_b$, we find a Hamilton path from $u_0$ to $u_b$ by removing the edges $u_0v_0$, $v_{2a}v_{2a-b}$, $u_{-b}v_{-b}$, $u_{2a+b}v_{2a+b}$, $u_bu_{a+b}$, $u_{a-b}u_{2a-b}$, and adding the edges $v_0v_{-b}$, $v_{2a}v_{2a+b}$, $u_{2a-b}v_{2a-b}$, $u_{-b}u_{a-b}$, $u_{a+b}u_{2a+b}$.
If the vertices  $v_{2a+b}$, $u_{2a+b}$ belong to the subpath $v_a\,P\,v_b$, then we find a Hamilton path from $v_{a-b}$ to $v_{2a-b}$ yielding a
Hamilton path from $v_0$ to $v_a$ by removing the edges $v_0v_b$, $v_{2a}v_{2a-b}$, $u_{-b}v_{-b}$, $u_bu_{a+b}$, $u_{a-b}v_{a-b}$, 
$u_{2a+b}v_{2a+b}$, and adding the edges $u_bv_b$, $v_0v_{-b}$, $v_{2a}v_{2a+b}$, $u_{-b}u_{a-b}$, $u_{a+b}u_{2a+b}$.
\medskip

\textbf{Case (b).} Assume that both edges $u_{-b}v_{-b}$ and $u_{2a+b}v_{2a+b}$ belong to the subpath $v_a\,P\,v_b$.
  If the vertices $v_{2a+b}$, $u_{2a+b}$  precede $u_{-b}$, $v_{-b}$ in $v_a\,P\,v_0$, then we find a Hamilton path from $v_{a-b}$ to $v_{2a-b}$ yielding a Hamilton path from $v_0$ to $v_a$ by removing the edges $v_0v_b$, $v_{2a}v_{2a-b}$, $u_{-b}v_{-b}$, $u_bu_{a+b}$, $u_{a-b}v_{a-b}$, $u_{2a+b}v_{2a+b}$, and adding the edges $u_bv_b$, $v_0v_{-b}$, $v_{2a}v_{2a+b}$, $u_{-b}u_{a-b}$, $u_{a+b}u_{2a+b}$. 
If the vertices $v_{2a+b}$, $u_{2a+b}$ follow $u_{-b}$, $v_{-b}$ in $v_a\,P\,v_0$, then we find a Hamilton path from $u_0$ to $u_b$ by
removing the edges $u_0v_0$, $v_{2a}v_{2a-b}$, $u_bu_{a+b}$, $u_{a-b}u_{2a-b}$, $u_{-b}v_{-b}$, $u_{2a+b}v_{2a+b}$, and adding the edges $v_0v_{-b}$, $u_{a+b}u_{2a+b}$, $v_{2a}v_{2a+b}$, $u_{2a-b}v_{2a-b}$, $u_{-b}u_{a-b}$.  
\medskip

\textbf{Case (c).} Assume that the edge $u_{-b}v_{-b}$  belongs to the subpath $v_a\,P\,v_b$ and the edge $v_{2a+b}$, $u_{2a+b}$ belongs to the subpath  $u_a\,P\,u_b$.
We will use the edges $u_{a+2b}v_{a+2b}$ and $u_{-a}v_{-a}$ to find the required Hamilton path. Note that these edges belong to the cycle $C$. If the vertices $u_{a+b}$, $v_{a+b}$, $u_{a+2b}$, $v_{a+2b}$ occur in this order in $C$, we find a Hamilton path from $u_{a+b}$ to $u_{a+2b}$ yielding a Hamilton path from $u_0$ to $u_b$ by removing the
edges $u_xv_x$ with $x\in\{a+b, a+2b\}$, and adding the edge $v_{a+b}v_{a+2b}$. Analogously, if the vertices $u_0$, $v_0$, $u_{-a}$, $v_{-a}$ occur in this order in $C$, we find a Hamilton path from $v_0$ to $v_{-a}$ yielding a Hamilton path from $v_0$ to $v_a$. Therefore, in the following, we assume that the vertices occur in the orders $u_{a+b}$, $v_{a+b}$, $v_{a+2b}$, $u_{a+2b}$, and $u_0$, $v_0$, $v_{-a}$, $u_{-a}$. Now we distinguish several cases with respect to the position of the edges $u_{-a}v_{-a}$ and $u_{a+2b}v_{a+2b}$ in $C$.
\medskip 

\textbf{First we consider the position of the edge $u_{-a}v_{-a}$ in $C$}, which could be in the path $u_a\,P\,u_b$ or in the path $v_a\,P\,v_b$.
If the edge $u_{-a}v_{-a}$ is in the path $u_a\,P\,u_b$ and precedes the edge $u_{2a+b}v_{2a+b}$, then we find a Hamilton path from $v_0$ to $v_{-a}$ providing a Hamilton path from $v_0$ to $v_a$ by removing the edges $u_xv_x$ with $x\in\{0, -a, -b, 2a+b\}$, $v_0v_b$, $v_{2a}v_{2a-b}$, $u_bu_{a+b}$, $u_{a-b}u_{2a-b}$,
and adding the edges $v_0v_{-b}$, $u_0u_{-a}$, $v_{2a}v_{2a+b}$, $u_{2a-b}v_{2a-b}$, $u_{a+b}u_{2a+b}$, $u_bv_b$, $u_{-b}u_{a-b}$. The same holds if the edge $u_{-a}v_{-a}$ belongs to the path $v_a\,P\,v_b$ and follows the edge $u_{-b}v_{-b}$. We find a Hamilton path from $v_0$ to $v_{-a}$, yielding a Hamilton path from $v_0$ to $v_a$, even if the path $(v_{-a+b}, v_{-a}, u_{-a})$ is in $v_a\,P\,v_b$ and precedes $u_{-b}v_{-b}$ by removing the edges $v_0v_b$, $u_bu_{-a+b}$, $v_{-a}v_{-a+b}$, and adding the edges $u_bv_b$, $u_{-a+b}v_{-a+b}$.
It remains to consider the case where $(v_{-a-b}, v_{-a}, u_{-a})$ is in $v_a\,P\,v_b$ and precedes $u_{-b}v_{-b}$,  and the case where $u_{-a}v_{-a}$ is in $u_a\,P\,u_b$ and follows the edge $u_{2a+b}v_{2a+b}$. 
\medskip

\textbf{Next we consider the position of the edge $u_{a+2b}v_{a+2b}$ on $C$}: it could be in $u_a\,P\,u_b$ or in $v_a\,P\,v_b$.
If the path $(v_{a+2b}, u_{a+2b}, u_{2b})$ is in $u_a\,P\,u_b$, then we find a Hamilton path from $u_{a+b}$ to $u_{a+2b}$ regardless of the position of the edge $u_{2a+b}v_{2a+b}$ with respect to $(v_{a+2b}, u_{a+2b}, u_{2b})$ by removing the edges $v_bv_{2b}$, $u_{2b}u_{a+2b}$,
$u_bu_{a+b}$, and adding the edges $u_bv_b$, $u_{2b}v_{2b}$; the path yields a Hamilton path from $u_0$ to $u_b$.
We find a Hamilton path from $u_{a+b}$ to $u_{a+2b}$ even if the path $(v_{a+2b}, u_{a+2b}, u_{2a+2b})$ is in $u_a\,P\,u_b$ and follows the edge $u_{2a+b}v_{2a+b}$ by removing the edges $v_0v_b$, $v_{2a}v_{2a-b}$, $u_{2a+b}v_{2a+b}$, $u_{a+2b}v_{a+2b}$, $u_bu_{a+b}$, $u_{a+b}v_{a+b}$, $u_{a-b}u_{2a-b}$, $u_{-b}v_{-b}$, and adding the edges  $u_bv_b$, $v_0v_{-b}$, $v_{2a}v_{2a+b}$, $u_{2a-b}v_{2a-b}$, $u_{a+b}u_{2a+b}$, $v_{a+b}v_{a+2b}$, $u_{-b}u_{a-b}$. The case where the path $(v_{a+2b}, u_{a+2b}, u_{2a+2b})$ is in $u_a\,P\,u_b$ and precedes the edge $u_{2a+b}v_{2a+b}$ will be considered later.

Now we deal with the edge $u_{a+2b}v_{a+2b}$ in $v_a\,P\,v_b$. A Hamilton path from $u_{a+b}$ to $u_{a+2b}$ can also be found even if the edge $u_{a+2b}v_{a+2b}$ in $v_a\,P\,v_b$ and precedes the edge $u_{-b}v_{-b}$ by removing the edges $u_xv_x$ with $x\in\{2a+b, a+b, a+2b, -b\}$, $v_0v_b$, $v_{2a}v_{2a-b}$, $u_bu_{a+b}$, $u_{a-b}u_{2a-b}$, and adding the edges $v_0v_{-b}$, $v_{2a}v_{2a+b}$, $u_{2a-b}v_{2a-b}$, $u_{a+b}u_{2a+b}$, $u_bv_b$, $v_{a+b}v_{a+2b}$, $u_{-b}u_{a-b}$.
We find a Hamilton path from $u_{a+b}$ to $u_{a+2b}$, or from $u_{-a+b}$ to $u_{-a+2b}$, yielding a Hamilton path from $u_0$ to $u_b$,
even if $u_{a+2b}v_{a+2b}$ follows $u_{-b}v_{-b}$ in $v_a\,P\,v_b$, and the path $(v_{2b}, u_{2b}, u_{a+2b}, v_{a+2b})$ is not contained in $C$.
More specifically, we have a Hamilton path from $u_{a+b}$ to $u_{a+2b}$ when the path $(u_{2a+2b}, u_{a+2b}, v_{a+2b})$ is in $v_a\,P\,v_b$:
remove the edges $v_0v_b$, $v_{2a}v_{2a-b}$, $v_{2a+b}v_{2a+2b}$, $u_bu_{a+b}$, $u_{a-b}u_{2a-b}$, $u_{-b}v_{-b}$, $u_{a+2b}u_{2a+2b}$, and add the edges $v_0v_{-b}$, $v_{2a}v_{2a+b}$, $u_{2a-b}v_{2a-b}$, $u_{2a+2b}v_{2a+2b}$, $u_{-b}u_{a-b}$, $u_bv_b$;
we have a Hamilton path from $u_{-a+b}$ to $u_{-a+2b}$ when the path $(u_{-a+2b}, u_{2b}, u_{a+2b}, v_{a+2b})$ is in $v_a\,P\,v_b$: 
remove the edges $v_bv_{2b}$, $u_{b}u_{-a+b}$, $u_{2b}u_{-a+2b}$, and add the edges $u_bv_b$, $u_{2b}v_{2b}$. It remains to consider the case where the path $(v_{2b}, u_{2b}, u_{a+2b}, v_{a+2b})$ is in $C$.

\textbf{By intersecting the open cases mentioned above}, the following cases remain to be studied:
(1) the path $(v_{a+2b}, u_{a+2b}, u_{2a+2b})$ precedes the edge $u_{2a+b}v_{2a+b}$ in $u_a\,P\,u_b$ and the edge $u_{-a}v_{-a}$
follows  the edge $u_{2a+b}v_{2a+b}$ in $u_a\,P\,u_b$, (2) the path $(v_{a+2b}, u_{a+2b}, u_{2a+2b})$ precedes the edge $u_{2a+b}v_{2a+b}$ in $u_a\,P\,u_b$ and $(v_{-a-b}, v_{-a}, u_{-a})$ precedes $u_{-b}v_{-b}$ in $v_a\,P\,v_b$, (3) the path $(v_{2b}, u_{2b}, u_{a+2b}, v_{a+2b})$ is in $C$ and the edge $u_{-a}v_{-a}$ follows  the edge $u_{2a+b}v_{2a+b}$ in $u_a\,P\,u_b$, (4) the path $(v_{2b}, u_{2b}, u_{a+2b}, v_{a+2b})$ is in $C$ and $(v_{-a-b}, v_{-a}, u_{-a})$ precedes $u_{-b}v_{-b}$ in $v_a\,P\,v_b$.

\textbf{Case (1).} Assume that the path  $(v_{a+2b}, u_{a+2b}, u_{2a+2b})$ precedes the edge $u_{2a+b}v_{2a+b}$ in $u_a\,P\,u_b$ and 
the edge $u_{-a}v_{-a}$ follows  the edge $u_{2a+b}v_{2a+b}$ in $u_a\,P\,u_b$. If the path $(u_{-a}, v_{-a}, v_{-a+b}, u_{-a+b})$ is in $C$,  then the path $(v_{2b-a},u_{2b-a},u_{2b},v_{2b})$ is in $C$, 
so the vertices $u_{-a+2b}$, $v_{-a+2b}$, $u_{-a+b}$, $v_{-a+b}$ occur in this order in $C$. Hence we can find a Hamilton path from
$u_{-a+2b}$ to $u_{-a+b}$ yielding a Hamilton path from $u_0$ to $u_b$: remove the edges $u_xv_x$ with $x\in\{-a+2b, -a+b\}$, and
add the edge $v_{-a+b}v_{-a+2b}$. If the path $(u_{-a}, v_{-a}, v_{-a+b})$ is in $C$ but the vertices $v_{-a+b},u_{-a+b}$ are not adjacent in $C$, then the path $(u_{-a}, v_{-a}, v_{-a+b}, v_{-a+2b})$ is in $C$, and we find a Hamilton path from $v_{2b}$ to $v_{-a+2b}$ yielding a Hamilton path from $v_0$ to $v_a$: remove the edges $v_bv_{2b}$, $v_{-a+b}v_{-a+2b}$, $u_bu_{-a+b}$, and add the edges
$u_bv_b$, $u_{-a+b}v_{-a+b}$. If $(u_{-a}, v_{-a}, v_{-a-b})$ is in $C$, then $u_{-a+b}v_{-a+b}$ is in $C$, and we can find
a Hamilton path from $v_0$ to $v_{-a}$ yielding a Hamilton path from $v_0$ to $v_a$: remove the edges $v_0v_b$,
$v_{2a}v_{2a-b}$, $u_{2a+b}v_{2a+b}$, $v_{-a}v_{-a-b}$, $u_bu_{a+b}$, $u_{a-b}u_{2a-b}$, $u_{-b}u_{-a-b}$, and add the edges
$u_bv_b$, $v_{2a}v_{2a+b}$, $u_{2a-b}v_{2a-b}$, $u_{a+b}u_{2a+b}$, $u_{-a-b}v_{-a-b}$, $u_{-b}u_{a-b}$.
The constructions for $(u_{-a}, v_{-a}, v_{-a+b}, u_{-a+b})$ not belonging to $C$ can also be repeated for the case where $(v_{2b}, u_{2b}, u_{a+2b}, v_{a+2b})$ is in $C$ and the edge $u_{-a}v_{-a}$ follows the edge $u_{2a+b}v_{2a+b}$ in $u_a\,P\,u_b$.

\textbf{Case (2).} Assume that the path $(v_{a+2b}, u_{a+2b}, u_{2a+2b})$ precedes the edge $u_{2a+b}v_{2a+b}$ in $u_a\,P\,u_b$ and 
the path $(v_{-a-b}, v_{-a}, u_{-a})$ precedes the edge $u_{-b}v_{-b}$ in $v_a\,P,v_b$. Notice that the edge $u_{-a+b}v_{-a+b}$ is in $C$. We find a Hamilton path from $u_{-a+b}$ to $u_{-a+2b}$ yielding a Hamilton path from $u_0$ to $u_b$ as follows: remove the edges 
$u_xv_x$ with $x\in\{a+2b, -a+b, -a\}$, $u_0u_a$, $v_av_{a+b}$, $u_{2b}u_{-a+2b}$, and add the edges
$u_0u_{-a}$, $u_av_a$, $v_{a+b}v_{a+2b}$, $u_{2b}u_{a+2b}$, $v_{-a}v_{-a+b}$.

\textbf{Cases (3)--(4).} Assume that the path $(v_{2b}, u_{2b}, u_{a+2b}, v_{a+2b})$ is in $C$, and the path $(u_{-a}, v_{-a}, v_{-a+b}, u_{-a+b})$ follows the edge $u_{2a+b}v_{2a+b}$ in $u_a\,P\,u_b$, or the path $(v_{-a-b}, v_{-a}, u_{-a})$ precedes the edge $u_{-b}v_{-b}$ in $v_a\,P,v_b$. Note that the cycle $C$ contains the edge $u_{-a+2b}v_{-a+2b}$ in both cases.
In the former case, we can assume that the vertices $v_{2b}$, $u_{2b}$, $u_{-a+2b}$, $v_{-a+2b}$ occur in this order in $C$,
otherwise we find a Hamilton path from $v_{2b}$ to $v_{-a+2b}$ providing a Hamilton path from $v_0$ to $v_a$ 
(remove the edges $u_xv_x$ with $x\in\{2b, -a+2b\}$, and add the edge $u_{2b}u_{-a+2b}$). Then the vertices 
$v_{-a+2b}$, $u_{-a+2b}$, $v_{-a+b}$, $u_{-a+b}$ occur in this order in $C$, and we can find a Hamilton path from $u_{-a+2b}$ to
$u_{-a+b}$ yielding a Hamilton path from $u_0$ to $u_b$: remove the edges $u_xv_x$ with $x\in\{-a+b, -a+2b\}$, 
and add the edge $v_{-a+b}v_{-a+2b}$. In the latter case, the subpath $(u_{-a+2b}, v_{-a+2b}, v_{-a+b}, u_{-a+b})$ is in $C$, and consequently the vertices 
$u_{2b}$, $v_{2b}$, $u_{-a+2b}$, $v_{-a+2b}$ occur in this order in $C$. We can thus remove the edges $u_xv_x$ with $x\in\{2b, -a+2b\}$,
and add the edge $u_{2b}u_{-a+2b}$. We find a Hamilton path from $v_{2b}$ to $v_{-a+2b}$ yielding a Hamilton path from $v_0$ to $v_a$.
\end{proof}

\end{document}